\documentclass[12pt]{article}
\usepackage[letterpaper, left=.5in, right=.5in, top=.75in, bottom=1in ]{geometry}

\usepackage[utf8]{inputenc}
\usepackage{color}
\usepackage{amsmath}
\usepackage{amsthm}
\usepackage{amssymb}
\usepackage{array,epsfig}
\usepackage{dsfont}
\usepackage{mathrsfs}
\usepackage{mathtools}
\usepackage{tikz}
\usepackage[numbers]{natbib}
\usepackage{graphicx}
\usepackage{verbatim}

\newtheorem{thm}{Theorem}[section]

\newtheorem{lem}[thm]{Lemma}

\newtheorem{prop}[thm]{Proposition}

\theoremstyle{definition}

\newtheorem{defn}[thm]{Definition}

\newtheorem{rmk}[thm]{Remark}

\newcommand{\ra}{\rightarrow}

\newcommand{\bb}[1]{\mathbb{#1}}

\newcommand{\Q}{\bb{Q}}
\newcommand{\R}{\bb{R}}

\newcommand{\PP}{\bb{P}}
\newcommand{\E}{\bb{E}}

\newcommand{\calM}{\mathcal{M}}

\newcommand{\calW}{\mathcal{W}}

\newcommand{\bPm}{\bar\PP_\mu}

\newcommand{\hPm}{\hat\PP_\mu}

\newcommand{\ttkd}{t \wedge \tau_k^\delta}
\newcommand{\tskd}{t \wedge \sigma_k^\delta}

\newcommand{\interior}[1]{%
	{\kern0pt#1}^{\mathrm{o}}%
}

\title{Reflected Brownian Motion with Drift in a Wedge}
\date{}
\author{Peter Lakner, Ziran Liu, Josh Reed}

\usepackage{natbib}

\usepackage{graphicx}

\usepackage{verbatim}

\begin{document}
	
	\maketitle
	
\centerline{Department of Technology, Operations, and Statistics, Stern School of Business}
\centerline{New York University}
\centerline{pl3@stern.nyu.edu, zliu@stern.nyu.edu, jr180@stern.nyu.edu }

\abstract{We study reflecting Brownian motion with drift constrained to a wedge in the plane. Our first set of results provide necessary and sufficient conditions for existence and uniqueness of a solution to the corresponding submartingale problem with drift, and show that its solution possesses the Markov and Feller properties. Next, we study a version of the problem with absorption at the vertex of the wedge. In this case, we provide a condition for existence and uniqueness of a solution to the problem and some results on the probability of the vertex being reached.}

\section{Introduction}
\label{submartingale problem}
	
	
	
In this paper, we study 2-dimensional Brownian motion with constant drift $\mu \in \mathbb{R}^2$ constrained to a wedge $S$ in $\mathbb{R}^2$. This process may also be referred to as reflected Brownian motion (RBM) with drift in a wedge, and we denote the process itself by $Z$. For concreteness, we define the wedge in polar coordinates by $\{r \geq 0, 0 \leq \theta \leq \xi \}$ for some $0 < \xi < 2 \pi$. Loosely speaking, the behavior of $Z$ may be characterized as follows. In the interior of $S$, $Z$ behaves as a 2-dimensional Brownian motion. On the other hand, the behavior of $Z$ on the boundary of $S$ is characterized by two reflection angles $\theta_1$ and $\theta_2$, depending upon whether the lower boundary $\partial S_1$ or upper boundary $\partial S_2$ has been reached. Both $0 < \theta_1, \theta_2, < \pi$ and the angles are measured from their inward facing normals,  with positive angles  corresponding to reflection toward the vertex of the wedge and negative angles away.   See Figure 
below for an illustration.
	  
\begin{figure}[h!]
	\centering
	\includegraphics[scale=0.38]{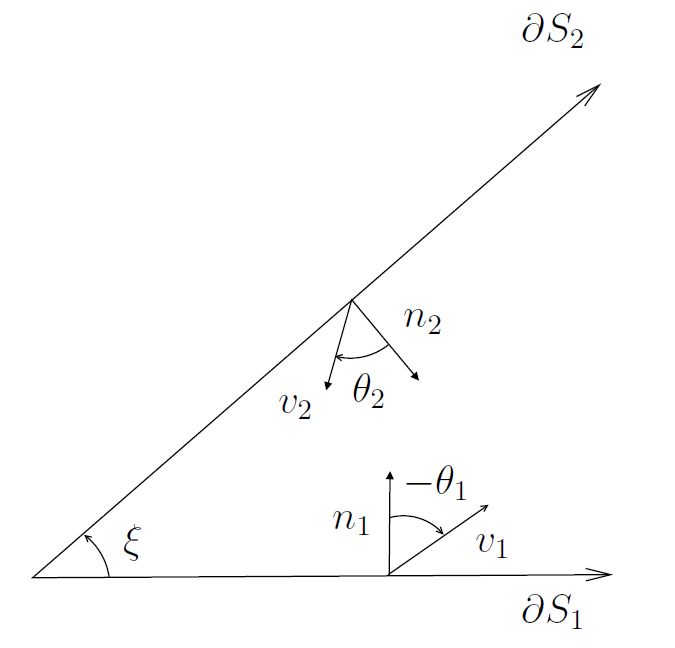}
	\caption{RBM in a Wedge}
	\label{fig:RBM in a Wedge}
\end{figure}

	One way to define RBM in a wedge is using a sample-path approach \cite{dai1992reflected,dai1996existence,harrison1981reflected,taylor1993existence} where $Z$ is explicitly characterized as the sum of a $2$-dimensional Brownian motion on an arbitrary probability space \cite{karatzas2012brownian,le2016brownian,revuz2013continuous} and a constraining or pushing process which satisfies the specifications related to the directions of reflection given above.  This sample-path approach works with or without a drift for some but not all parameter regimes of $(\xi, \theta_1,\theta_2)$. It tends  not to work in regimes where $Z$ is known not to be a semi-martingale \cite{williams1985reflected} and the pushing process has infinite variation. Recent progress in this direction has however been made \cite{kang2010dirichlet,ramanan2006reflected}. 
	
	A more probabilistic approach to defining RBM in a wedge was given by Varadhan and Williams \cite{varadhan1985brownian}. In this case, $Z$ is defined as the solution to a submartingale problem. This approach yields existence and uniqueness results for all parameter regimes but at several points the proofs of \cite{varadhan1985brownian}  rely heavily on the assumption that $Z$ behaves as a standard Brownian motion inside of $S$.  This is not an issue for parameter regimes where the sample-path approach described above may be applied because it is amenable to Brownian motions with drift, and the recent paper \cite{kang2017submartingale} demonstrates equivalence between the sample-path and the submartingale approach in such settings. On other hand, in parameter regimes where the sample-path approach cannot be applied, extending the results of \cite{varadhan1985brownian} in the direction of allowing $Z$ to behave as a Brownian motion with drift in the interior of $S$ remains an open problem. In this paper, we resolve this open problem. 
	
	Our primary motivation comes from queueing theory where semi-martingale RBM with drift has long been known to serve as the weak limit of both the properly scaled queue length \cite{chen2002brownian,harrison1978diffusion,harrison1981distribution,reiman1984open} and workload \cite{bramson2001heavy,chen1996diffusion,peterson1991heavy,williams1998diffusion} processes of different queueing systems in heavy-traffic. In such queueing settings, the drift term arises as the result of an imbalance between the input and output processes to the system. The limiting RBM in these cases is often defined using the sample-path approach via the conventional Skorokhod map \cite{harrison1981reflected,skorokhod1961stochastic,whitt2002stochastic}. More recently, using the extended Skorokhod map \cite{ramanan2006reflected}, RBM with drift which is not a semi-martingale has been proven \cite{ramanan2003fluid} to be the weak limit of the properly scaled unfinished work process of the generalized processor sharing model in heavy traffic. In this example, the sample-path approach is still employed to define the limiting process with the help of the extended Skorokhod map \cite{ramanan2006reflected}. We conjecture however that there exist other applied queueing settings where the limiting heavy-traffic process is an RBM with drift which is not a semi-martingale and cannot be rigorously defined via the sample-path approach. One of these settings is the coupled processor model \cite{cohen2000boundary,fayolle1979two}. In such situations, before proving any limit theorems, it is necessary to first establish the existence of RBM with drift through other means such as the submartingale problem.
	
	The remainder of the paper is organized as follows. Our main results may be found in Section \ref{main results}. In Section 	\ref{sec:sub:with:drift:results},  we provide necessary and sufficient conditions for the existence and uniqueness of the solution to the submartingale problem with drift (see Definition \ref{defn of the submartingale problem}), and show that its solution possess the strong Markov property and three versions of the Feller property. Next, in Section \ref{sec:absorbed:with:drift:results}, we study the submartingale problem with drift absorbed at the vertex of the wedge (see Definition \ref{definition of absorbed process}). We provide results on the  existence and uniqueness of the solution to this problem and results on the probability of the absorbed process with drift reaching the vertex of the wedge. Sections \ref{Existence} through \ref{sec:hitting} contain the proofs of our main results.

	\section{Main Results}\label{main results}
	
	Before stating our main results, we first set up some notation. Let $C_S = C(\R_+, S)$ and, for each $t \ge 0$, let  $Z(t): C_S \ra S$ denote the coordinate map $Z(t)(\omega) = \omega(t)$ for $\omega \in C_S$. Also, let $Z = \{Z(t), t\ge 0\}$ denote the coordinate mapping process on $C_S$. Let $\calM_t  = \sigma(Z(s), 0\le s\le t)$ be the underlying natural filtration with terminal $\sigma$-algebra $\calM = \sigma (Z(s), s\ge 0)$.  For each $n \ge 1$ and domain $\Omega \subseteq \R^2$, we denote by $C^n_b(\Omega)$ the set of $n$ times bounded continuously differentiable functions on $\Omega$. We assume that the wedge $S$ is positioned so that one side of it is the positive horizontal half line, and the angle of the wedge is $\xi$. We define $\partial S_1$ and $\partial S_2$ as the two sides of the wedge so that neither includes the vertex, i.e., $\partial S_1=\{(x,0):\ x>0\}$ and $\partial S_2=  \{r(\cos\xi,\sin\xi):\ r>0\}$.  Next (see Figure \ref{fig:RBM in a Wedge}), we denote by $v_1$ and $v_2$ the reflection directions on the boundaries $\partial S_1$ and $\partial S_2$, respectively. For convenience, we assume that each $v_i$ is normalized such that $v_i \cdot n_i = 1$, where $n_i$ is the inward facing normal vector on $\partial S_i$ for $i=1,2$. Finally, for $i=1,2$, we set the directional derivative operator  $D_i = v_i \cdot \nabla$, with $\nabla$ being the gradient operator, the dot is the inner product,  and denote by $\Delta$ the Laplacian operator.

	\subsection{The Submartingale Problem with Drift}
	\label{sec:sub:with:drift:results}
	
	\begin{defn}[Submartingale Problem with Drift]\label{defn of the submartingale problem}
		
		A family of probability measures $\{\PP_{\mu}^z, z \in S\}$ on $(C_S, \mathcal{M})$ is said to solve the submartingale problem with drift $\mu \in \R^2$ if for each $z \in S$, the following three conditions hold,
		
		\begin{itemize}
			\item[1.] $\PP_{\mu}^z(Z(0)=z) = 1$;
			\item[2.] For each $f \in C_b^2(S)$, the process $$\bigg\{f(Z(t)) - \int_0^t \mu \cdot \nabla f(Z(s))ds -\frac{1}{2} \int_0^t \Delta f(Z(s))ds, t \ge 0\bigg\}$$ is a submartingale on $(C_S, \calM, \calM_t, \PP_\mu^z)$ whenever $f$ is constant in a neighbourhood of the origin and satisfies $D_i f \ge 0$ on $\partial S_i$ for $i = 1, 2$;
			\item[3.] $\E_\mu^z\bigg[\displaystyle \int_0^\infty \mathds{1}_{\{Z(t) = 0\}} dt \bigg] = 0$.
		\end{itemize}
	\end{defn}

The above definition bears a relationship to  the extended Skorokhod problem (ESP)  developed in \cite{ramanan2006reflected}. We shall recall the definition of the ESP below. Let $d(\cdot)$ be a set-valued map from $\partial S$, the boundary of $S$, to the class of subsets of $\R^2$ satisfying the following two conditions:
	\begin{itemize}
		\item[(d1)] for any $x \in \partial S$, the image $d(x)$ is a non-empty closed convex cone in $\R^2$ with the vertex being the origin;
		\item[(d2)] the graph $\{(x, d(x)); \ x \in \partial S\}$ is closed.
	\end{itemize}
	
	For convenience, we  extend the definition of $d(\cdot)$ to $S$ by setting $d(x) = \{0\}$ for all $x \in S^\circ$, where $S^\circ$ is the interior of $S$. For a set $A\subset {\mathbb R}^2$, let ${\rm co}(A)$ be the closed convex cone generated by $A$.

\begin{defn}[Extended Skorokhod Problem (ESP)]\label{ESP}
		A pair of processes $(\phi, \eta) \in C_S\times C(\R_+,\R^2)$ is said to solve the ESP $(S, d(\cdot))$ for $\psi \in C(\R_+,\R^2)$ such that $\psi(0)\in S$ if $\phi(0) = \psi(0)$, and if for all $t \in \R_+$, the following properties hold,
		\begin{enumerate}
			\item $\phi(t) = \psi(t) + \eta(t)$; 
			\item $\phi(t) \in S$; \label{item:2:esp}
			\item For every $s \in [0,t]$;
			$$
			\eta(t)-\eta(s) \in \hbox{co}\Big[\cup_{u\in(s,t]}d(\phi(u))\Big].
			$$	
		\end{enumerate}
	\end{defn}
Item \ref{item:2:esp} in the above definition is redundant since we already required that $\phi\in C_S$, but we kept that item as it appears in the original definition in \cite{ramanan2006reflected}.
	
	Just like  in \cite{varadhan1985brownian}, let
	$$
	\alpha = \frac{\theta_1 + \theta_2}{\xi}.
	$$
	The quantity $\alpha$ plays a prominent role.\\
	
	\begin{thm}\label{Exsistence and Uniqueness of Submartingale Problem}
		If $\alpha <2$, then for each $\mu \in \mathbb{R}^2$  there exists a unique solution $\{\PP^z_\mu,z\in S\}$ to the submartingale problem with drift. In addition, the following statements hold:
	\begin{itemize}
			\item[1.] There exists a process $X$ defined on $(C_S, {\calM}, {\calM}_t)$ which, for each $z \in S$, is a  $2$-dimensional Brownian motion with drift $\mu$ started at $z$ under $\PP^z_\mu$;
			\item[2.] Setting $Y=Z-X$, the pair $(Z,Y)$ solves the ESP $(S,d(\cdot))$ for $X$, $\PP^z_\mu$-a.s., for each $z\in S$.
		\end{itemize}

	\end{thm}
	The above theorem establishes a decomposition
	\begin{equation} Z=X+Y,\label{decomposition:main}
	    \end{equation}
	such that for all $z\in S$ under $\PP_\mu^z$ the process $X$ is a standard Brownian motion with drift $\mu$ started at $z$. In the following two theorems we shall establish several  properties of the process $Y$ appearing in the above decomposition.
	In order to state the first of these two  results, we need the  definition of the strong $p$-variation of a function. Let $T>0$ arbitrary. We call an ordered set $(t_0,t_1,\dots,t_n)$ a partition of the interval $[0,T]$, if $0=t_0<t_1<\dots <t_n=T$, for an arbitrary $n\in{\mathbb N}_+$.  Let $\pi(T)$ denote the set of all partitions of the interval $[0,T]$. We define the mesh of a partition $\rho=(t_0,\dots,t_n)\in\pi(T)$ by setting
	$$\|\rho\|=\max\{t_i-t_{i-1}: i=1,\dots,n\}.   $$
	
	\begin{defn}
	        Let $T>0$ and $p>0$. The strong $p$-variation of a function $f:{\mathbb R}_+\mapsto {\mathbb R}^k$ on $[0,T]$ is defined by
	        $$V_p(f,[0,T])= \sup\left\{\sum_{t_i\in\rho,\ i\ge 1} \|f(t_i)-f(t_{i-1}\|^p: \rho\in\pi(T)\right\}.$$
	\end{defn}

	\begin{thm}\label{variation} Suppose that $1<\alpha<2$. Then for each $p>\alpha$ and $z\in S$,
	\begin{equation} \PP_\mu^z(V_p(Y,[0,T])<+\infty)=1,\quad T>0,\label{var_finite}\end{equation}
	and, for each $0<p\le \alpha$,
	\begin{equation} \PP_\mu^0(V_p(Y,[0,T])<+\infty)=0,\quad T>0.\label{var_infinite}\end{equation}
	    
	\end{thm}
	\vskip .1in
	
	A 2-dimensional continuous process $U$ defined on $(C_s,{\mathcal M}, {\mathcal M}_t, \PP^z_\mu)$ is said to be of {\it zero energy} if for each $T>0$ and each sequence of partitions $(\rho^m)\subset \pi(T)$ such that $\|\rho^m\|\to 0$ as $m\to\infty$ we have
	$$\sum_{i=1}^{n(m)}\|U(t_i^m) -U(t_{i-1}^m\|^2\buildrel{ \PP^z_\mu}\over\longrightarrow 0\ {\rm as}\ m\to\infty, $$
	where
	$\rho^m =(t_0^m,\dots, t_{n(m)}^m)$.
	A process $D$ on $(C_s,{\mathcal M}, {\mathcal M}_t, \PP^z_\mu)$ is said to be a {\it Dirichlet process} if it has a decomposition $D=M+U$, where $M$ is a local martingale on the same probability space, and $U$ is a continuous zero-energy process with $U(0)=0$.
	
	\begin{thm}\label{dirichlet} Let $1<\alpha<2$, and $z\in S$ arbitrary. Then the process $Y$ in  decomposition 
	\eqref{decomposition:main} is a zero-energy process, and $Z$ is a Dirichlet process.     
	    
	\end{thm}

	\begin{thm}\label{notsemimartingale} If $1\le \alpha<2$, then $Z$ is not a semimartingale on $(C_S, \calM, \calM_t, \PP_\mu^z)$, for any $z\in S$.

	\end{thm}

	\begin{thm}\label{nonexistence} If $\alpha\ge 2$, then for any $\mu\in {\mathbb R}^2$ there is no solution to the submartingale problem with drift.

		\end{thm}
	
	Let $\{\PP_{\mu}^z, z \in S\}$ be the solution to the submartingale problem for some $\mu \in \mathbb {R}^2$. We say that 
	$\{\PP_{\mu}^z, z \in S\}$ possesses the strong Markov property if for each stopping time $\tau$ and $z \in S$, and each bounded $\mathcal{M}$-measurable function $h:C_S \ra \mathbb{R}$ we have that
	\begin{eqnarray}
		\E_\mu^z[\mathds{1}_{\{\tau < \infty\}}f(\omega(\cdot+\tau))|\calM_\tau] = \mathds{1}_{\{\tau < \infty\}} \E_\mu^{\omega(\tau)}[f(w(\cdot))], \,\,\,\,\,\,\, \PP_\mu^z\text{-a.s.} 
	\end{eqnarray}
	
	\begin{thm}\label{Strong Markov for Submartingale Problem}
		If $\alpha <2$, then for each $\mu \in \mathbb{R}^2$  the solution to the submartingale problem with drift has the strong Markov property.
	\end{thm}
	


	
	The last subject of this subsection is the Feller property of $\{\PP_{\mu}^z, z\in S\}$. There are various, slightly differing definitions of the Feller property available  in the literature. For clarity we list below three definitions.\\
	
	\noindent 1. We say that $\{\PP_{\mu}^z,  z \in S\}$ has the Feller property if for any $\{z_n, n \ge 1\} \subset S$ converging to $z \in S$, $\PP_\mu^{z_n} \Rightarrow \PP_\mu^z$ as $n \ra \infty$ (see Varadhan and Williams \cite{varadhan1985brownian}).\\

	\noindent 2.  Let $\hat{C}(S)$ is the set of continuous functions on $S$ vanishing at infinity. We say that $\{\PP_{\mu}^z,  z \in S\}$  has the $\hat{C}(S)$-Feller property if 
	for any $f \in \hat{C}(S)$,  and $t \ge 0$, the function  $z\mapsto\E_\mu^z[f(Z_t)]$ is also in $\hat{C}(S)$.\\
	
	\noindent 3. We say that $\{\PP_{\mu}^z,  z \in S\}$  has the ${C}_b(S)$-Feller property if 
	for any $f \in {C}_b(S)$,  and $t \ge 0$, the function  $z\mapsto\E_\mu^z[f(Z_t)]$ is also in ${C}_b(S)$.\\
	
	\begin{rmk}
		The Feller property obviously implies the $C_b(S)$-Feller property. The $\hat{C}(S)$-Feller property implies the $C_b(S)$-Feller, but the converse is not true (cf. Theorems 1.9 and 1.10 in \cite{bottcher2013levy}).\\
	\end{rmk}

	\begin{thm}\label{VW's Feller}
		If $\alpha <2$, then the solution to the submartingale problem for each $\mu \in \mathbb{R}^2$ has the Feller property.
	\end{thm}

	\begin{thm}\label{C-hat Feller}
		If $\alpha <2$, then the solution to the submartingale problem for each $\mu \in \mathbb{R}^2$ has the $\hat C(S)$-Feller property.
	\end{thm}

	We note that for the $\mu=0$ case the Feller property is known (\cite{varadhan1985brownian}, Theorem 3.13). However, the $\hat C(S)$-Feller property is new even in the $\mu=0$ case. 
	
	

	%

	\subsection{The Absorbed Process}
		\label{sec:absorbed:with:drift:results}
	
	Let
	$$
	\tau_0 = \inf \{t \ge 0: Z(t) = 0\}
	$$
	be the stopping time with respect to $\{\mathcal{M}_t, t \geq 0\}$ representing the first time that $Z$ reaches the vertex of the wedge.  Results in this subsection concern the RBM in a wedge up until $\tau_0$. Results of this type were provided in  \cite{varadhan1985brownian} for the driftless case and we extend many of them to the case of a constant drift $\mu$. We begin with the following definition.
	
	\begin{defn}[The Absorbed Process Problem]\label{definition of absorbed process}
		
		A family of probability measures $\{\PP_{\mu}^{z,0}, z \in S\}$ on $(C_S, \mathcal{M})$ is said to solve the absorbed process problem with drift $\mu \in \R^2$ if for each $z \in S$, the following three conditions hold,

		\begin{itemize}
			\item[1.] $\PP_{\mu}^{z,0}(Z(0)=z) = 1$;
			\item[2.] The process 
			$$
			\bigg\{f(Z({t\wedge \tau_0})) - \int_0^{t\wedge \tau_0} \mu \cdot \nabla f(Z(s))ds -\frac{1}{2} \int_0^{t\wedge \tau_0} \Delta f(Z(s))ds, t \ge 0\bigg\}
			$$
			is a submartingale on $(C_S, \calM, \calM_{t}, \PP_\mu^{z,0})$, for each $f \in C_b^2(S)$ such that $D_i f \ge 0$ on $\partial S_i$ for $i = 1, 2$;
			\item[3.] $\PP_{\mu}^{z,0}(Z(t)=0, \forall t \ge \tau_0) = 1$.
		\end{itemize}
		
	\end{defn}

	\begin{thm}\label{Existence and Uniqueness of Absorbed Process}
		For each $\mu \in \mathbb{R}^2$ and $\alpha\in{\mathbb R}$,  there exists a unique solution to the absorbed process problem. 
	\end{thm}
	
	The above theorem is  particularly interesting if $\alpha\ge 2$, since Theorem \ref{Exsistence and Uniqueness of Submartingale Problem} does not cover that case. The existence of a solution to the absorbed process problem easily follows from the existence of a solution to the submartingale problem whenever $\alpha < 2$. However, the uniqueness part of Theorem \ref{Existence and Uniqueness of Absorbed Process} does not follow in an obvious way from Theorem \ref{Exsistence and Uniqueness of Submartingale Problem} even in the $\alpha <2$ case. Our proof for Theorem \ref{Existence and Uniqueness of Absorbed Process} applies to all $\alpha\in{\mathbb R}$.\\

	Next we state a series of results on the hitting probability of the vertex for the absorbed process in the case of a constant drift. \\

	\begin{thm}\label{Hitting prob for alpha <= 0}
		If $\alpha \le 0$, then for each $\mu \in \R^2$ and $z \in S$, $\PP_\mu^{z,0}(\tau_0 = \infty) = 1$.
	\end{thm}

	The hitting probability of the vertex is more varied in the case of $\alpha \geq 1$, and before proceeding we must make some observations on the geometry of the wedge. For $n\ge 1$ and a set of vectors $\{a_1,\dots,a_n\}\subset\mathbb{R}^2$, let ${\rm co}(a_1,\dots,a_n)$ denote the closed convex cone generated by $\{a_1,\dots,a_n\}$.  We illustrate two relevant cases for $\alpha$ below.
	\vskip .2in
	
	\usetikzlibrary{arrows,quotes,angles}
	\usetikzlibrary{decorations.markings}
	\tikzstyle arrowstyle=[scale=1]
	\tikzstyle directed=[postaction={decorate,decoration={markings,
			mark=at position 1 with {\arrow[arrowstyle]{stealth}}}}]
	\tikzstyle reverse directed=[postaction={decorate,decoration={markings,
			mark=at position 1 with {\arrowreversed[arrowstyle]{stealth};}}}]
	\begin{tikzpicture}

		

		\draw
		(7.5,1.2)  coordinate (S) node[right] {$S$}
		(8.5,0) coordinate (a1) 
		-- (5,0) coordinate (b1) 
		-- (7.5,2.5) coordinate (c1) 
		pic["$\xi$", draw=black, <->, angle eccentricity=1.2, angle radius=1cm]
		{angle=a1--b1--c1};

		\coordinate (O) at (5,0) ;
		\coordinate (v1) at (4,1.2) ;
		\coordinate (v2) at (6,-1.2) ;
		
		\draw[directed] (O) -- (v1);
		\draw[directed] (O) -- (v2);
		\draw (4.3,0.8) node[right] {$v_1$};
		\draw (4.8,-0.8) node[right] {$v_2$};
		\draw (6,-2) node[right] {A};

		\draw
		(16.5,1.2)  coordinate (S) node[right] {$S$}
		(17.5,0) coordinate (a1) 
		-- (14,0) coordinate (b1) 
		-- (16.5,2.5) coordinate (c1) 
		pic["$\xi$", draw= black, <->, angle eccentricity=1.2, angle radius=1cm]
		{angle=a1--b1--c1};

		\coordinate (O) at (14,0) ;
		\coordinate (v1) at (14.5,1.2) ;
		\coordinate (v2) at (15.5,-0.8) ;
		
		\draw[directed] (O) -- (v1);
		\draw[directed] (O) -- (v2);
		\draw (14.3,0.8) node[left] {$-v_2$};
		\draw (13.8,-0.5) node[right] {$-v_1$};
		\draw (15,-2) node[right] {B};
	\end{tikzpicture}


	\noindent In the above diagrams, case A corresponds to $\alpha=1$, which occurs if and only if $\hbox{co}(v_1,v_2)$ is a line.  Case B corresponds to $\alpha>1$, which occurs  if and only if $\hbox{co}(-v_1,-v_2)$ contains $S$. In both cases, that is, whenever $\alpha\ge 1$ we have that ${\rm co}(v_1,v_2)\cap S =\{0\}.$ Note also that $\alpha \ge 1$ implies $\xi < \pi$.
	
	\begin{thm}\label{Hitting prob for alpha >= 1}
		If $\alpha \geq 1$, then 
		\begin{eqnarray}
			\PP_\mu^{z,0} (\tau_0<\infty) > 0 ~\mathrm{for~each}~ \mu \in \R^2 ~\mathrm{and}~ z \in S.  \label{finite}
		\end{eqnarray}
		Moreover, if in addition to the $\alpha\ge 1$ condition we also have that
		 \begin{equation}\hbox{co} (v_1,v_2,\mu) \cap S =\{0\},\label{condition}\end{equation}
		then for each $z \in S$, 
		\begin{equation}\PP_\mu^{z,0} (\tau_0<\infty)=1.\label{reaches}\end{equation}
	\end{thm}
	
	We note that in the case of $\alpha>1$ condition \eqref{condition} can be cast in an algebraic form. Let $R$ be the $2\times 2$ matrix such that its i-th column vector is $v_i$, for $i=1,2$. If $\alpha>1$ then condition \eqref{condition} is equivalent to the requirement  the vector $R^{-1}\mu$ has at least one non-negative component.
	
\begin{rmk} Theorem \ref{Hitting prob for alpha >= 1} leaves open the possibility that   $\PP^{z,0}_{\mu} (\tau_0<\infty) =1 $ whenever $\alpha \geq 1$. This however is not the case as the following counterexample shows. Let $\alpha \in \mathbb{R}$ be arbitrary and let the drift $\mu \in \mathbb{R}^2$ be given by $\mu = ||\mu||(\cos \eta, \sin \eta) \neq 0$, where  $\eta \in (0,\xi)$. Then, it is not hard to show using the proposition below that $\PP^{z,0}_{\mu} (\tau_0<\infty) < 1$ for each $z \in S\setminus\{0\}$. 
\end{rmk}

\begin{prop}\label{stays}Let $S$ be the 2-$d$ wedge defined above, let $S^0$ be the interior of $S$, let $B$ be a 2-$d$ standard Brownian motion with zero drift started at the origin under a probability measure $P$, and let  $\mu \in \mathbb{R}^2$ given by $\mu = ||\mu||(\cos \eta, \sin \eta) \neq 0$, where  $\eta \in (0,\xi)$. Then, if $0 < \xi < \pi$, for each $z \in \mathcal{S}^0$, 
\begin{eqnarray*}
P( z + B_t + \mu t \in \mathcal{S}^0, t \geq 0) > 0.
\end{eqnarray*}
\end{prop}

\noindent Using the proposition above and Theorem \ref{Hitting prob for alpha >= 1}, we may now deduce that if $\alpha \geq 1$, $\eta\in(0,\xi)$,  and $\mu = ||\mu||(\cos \eta, \sin \eta) \neq 0$, we obtain that $$\PP^{z,0}_{\mu} (\tau_0<\infty) \in (0,1)~\mathrm{for~each}~ z \in S\setminus\{0\}.$$ This implies that when $\alpha \geq 1$, hitting the vertex is no longer a $0$-$1$ event for certain values of $\mu$, which contrasts with the driftless result of \cite{varadhan1985brownian}.

\section{Proof of Theorems \ref{Exsistence and Uniqueness of Submartingale Problem}, \ref{variation}, \ref{dirichlet},  \ref{notsemimartingale}, and \ref{nonexistence}} \label{Existence}

	We  provide an extension of Theorems 2.4 and 2.8  in \cite{lakner2019roughness} to all $\alpha < 2$.  For $z \in \partial S_i, i = 1, 2,$ let $d(z) = \{\lambda v_i, \ \lambda \ge 0\}$, and set $d(0)=\R^2$. It is known that in the case of $\mu=0$ the submartingale problem has a unique solution whenever $\alpha<2$ (see \cite{varadhan1985brownian}). In accordance with our notation, that solution will be denoted by $\{\PP_0^z, z\in S\}$. 
	We then have the following.

	\begin{prop}\label{Theorem 2.4 in lakner2019roughness}
		Let $(C_S, {\calM}, {\calM}_t)$ and $Z$ be defined as in Section \ref{main results}. Then, if $\alpha < 2$,
		\begin{itemize}
			\item[1.] There exists a process $X$ defined on $(C_S, {\calM}, {\calM}_t)$ which, for each $z \in S$, is a  $2$-dimensional Brownian motion started at $z$ under $\PP^z_0$;
			\item[2.] Setting $Y=Z-X$, the pair $(Z,Y)$ solves the ESP $(S,d(\cdot))$ for $X$, $\PP^z_0$-a.s..
		\end{itemize}
	\end{prop}
	
	\begin{proof} Let $\alpha < 2$. Then, Condition 1 is immedate from Theorem 2.4 in \cite{lakner2019roughness}. It remains to show that for each $z \in S$, $\mathbb{P}^z_0$-a.s., $Z$ and $Y=Z-X$ together solve the ESP (see Definition \ref{ESP} as above) for $X$ with $d$ as defined immediately preceding the statement of the proposition. For $\alpha \in (1,2)$, this follows by Theorem 2.8 in \cite{lakner2019roughness}. We now claim that if $\alpha \le 1$, $(Z,Y)$ also solves the ESP $(S,d(\cdot))$ for $X$, $\PP^z$-a.s. For any two real numbers $0<s<t$, if $(s,t)$ belongs to a single excursion from the origin then by a similar proof to the one
		  in part 2 of Theorem 4.2 in \cite{lakner2019roughness}, one can conclude that item 3 in the definition of the ESP holds. If $(s,t)$ doesn't belong to one excursion, then item 3 is obviously satisfied by $d(0)={\mathbb R}^2$.
	\end{proof}
	
	
	
	We are now ready to prove the existence of a solution to the submartingale problem with drift, and some of the properties of the solution we create.
	
	\begin{prop}\label{existence theorem}
		If $\alpha <2$, then for each $\mu \in \R^2$ the submartingale problem with drift  has a solution $\{\PP_\mu^z, z\in S\}$, which satisfies the following properties. With $X$ and $Y$ defined in Proposition \ref{Theorem 2.4 in lakner2019roughness}, for every $z\in S$ the following hold:
	\begin{itemize}
	    \item Under $\PP^z_\mu$ the process $X$ is a standard Brownian motion with drift $\mu$ started at $z$; 
		\item The pair $(Z,Y)$ solves the ESP $(S,d(\cdot))$ for $X$, $\PP_\mu^z$-s.s. 
		\end{itemize}
	\end{prop}
	
	\begin{proof}
		
		
		Let $\alpha < 2$ and note that by {Proposition}  \ref{Theorem 2.4 in lakner2019roughness}  there exists a process $X$ defined on $(C_S, {\calM}, {\calM}_t)$ which, for each $z \in S$, is a  $2$-dimensional Brownian motion started at $z$ under $\PP^z_0$. Now
		let  $T \geq 0$ and for each $z\in S$, let $\PP_{\mu,T}^z$ be a probability measure on $(C_S, {\calM}, {\calM}_t)$ equivalent (mutually absolutely continuous) to $\PP^z_0$ such that under $\PP_{\mu,T}^z$, $X$ is a standard Brownian motion with drift $\mu$ up to time $T$, started at $z$. In other words, $\{X(t)-\mu t, t\le T\}$ is a standard (driftless) $\PP_{\mu,T}^z$-Brownian motion started at $z$. The measure $\PP_{\mu,T}^z$ is defined by 
		\begin{equation}
		\frac{d\PP_{\mu,T}^z} {d {\PP_0^z}} =\zeta(T), 
		\label{rn_derivative}\end{equation}
		where
		$\zeta(T)=\exp\{\mu\cdot(X(T)-X(0))-{1\over 2}\|\mu\|^2T\}$. 
		
		One can easily show that the family of probability measures $\{\PP_{\mu,T}^z,T\in[0,\infty)\}$ is consistent. That is, if $S<T$, then $\PP_{\mu,T}^z(A)=\PP_{\mu,S}^z(A)$, whenever $A\in{\calM_S}$.
		From \cite{parthasarathy2005probability}, Theorem 4.2 (page 143), it
		follows that there exists a single probability measure $\PP_{\mu}^z$ such that $\PP_{\mu}^z(A)=\PP_{\mu,T}^z(A)$ whenever $A\in {\cal M}_T$. 
		Since $\{X(t)-X(0)-\mu t, t\le T\}$ is a $\PP_{\mu,T}^z$-Brownian motion started at zero for every $T\in[0,\infty)$, it follows that $\{X(t)-X(0)-\mu t, t<\infty\}$ is also a  $\PP_{\mu}^z$-Brownian motion started at zero. Also,   $(Z,Y)$ solves the ESP $(S,d(\cdot))$ for $X$ under $\PP_{\mu}^z$-a.s., because by Proposition \ref{Theorem 2.4 in lakner2019roughness} it is true under $\PP_0^z$, and the measures $\PP_0^z$ and $\PP^z_\mu$ constrained to ${\mathcal M}_T$ are mutually absolutely continuous for every $T\in[0,\infty)$.
		Now let 
		$$W(t)=X(t)-X(0)-\mu t ,\quad t\in[0,\infty),$$
		and consider the definition of a weak solution to an SDER (see Definition 2.4 of Kang and Ramanan \cite{kang2017submartingale}).
		Clearly, the triplet $(C_S, {\cal M}, {\cal M}_t)$, $\PP_{\mu}^z$, $(Z,W)$ is a weak solution to the SDER with initial condition $z$ associated with $(S, d(\cdot))$, $b(\cdot)$ and $\sigma(\cdot)$, where $b(x)=\mu$ and $\sigma(x)=\hbox {Id}_{2\times 2}$. We note that the ``closed graph condition" (see Kang and Ramanan \cite{kang2017submartingale}, page 5) is satisfied. From Theorem 2 in \cite{kang2017submartingale}, it now follows that $\{\PP_{\mu}^z,z\in S\}$ solves the submartingale problem with drift $\mu$.
		
		
	\end{proof}


We shall use  the Lemmas \ref{Occupation:time}, \ref{w:delta}, \ref{construct the Brownian motion using the convergence}, \ref{Existence of BM under drift measure}, and  \ref{unique:main}
 for the proof of  both the uniqueness part of Theorem \ref{Exsistence and Uniqueness of Submartingale Problem}, and for the proof of Theorem \ref{nonexistence}. In these lemmas $\alpha$ may be an arbitrary real number. On the other hand, in these lemmas we start with a probability measure $P^z_\mu$ that satisfies conditions 1,2, and 3 of Definition \ref{defn of the submartingale problem}. This may be surprising, since Theorem \ref{nonexistence} states that such probability measure does not exists for $\alpha\ge 2$. However, for such $\alpha$'s we use these lemmas to derive a contradiction, thereby proving Theorem \ref{nonexistence}.
	
	\begin{lem}\label{Occupation:time}
		Suppose that $\{\PP_\mu^z, z\in S\}$ is a solution  to the submartingale problem with drift $\mu \in \R^2$. Then, for all $z \in S$, 
		
		\begin{equation}\E_\mu^z\bigg[\displaystyle \int_0^\infty \mathds{1}_{\{Z(t) \in \partial S\}} dt \bigg] = 0.\label{zero}\end{equation}
		
	\end{lem}

	\begin{proof}

	Let $z\in S$ be an arbitrary.	In this proof we shall use the Doob-Meyer decomposition for submartingales, which requires that the probability space is augmented. For this reason we denote by $(C_S,{\cal F}^z,({\cal F}_t^z))$ the augmentation of the space $(C_S,{\cal M},({\cal M}_t))$under $\PP_\mu^z$. For some technical details on the augmentation of probability spaces see Remark \ref{aug} in the Appendix. 
Condition $3$ of the submartingale problem gives
		$$
		\E_\mu^z\bigg[\displaystyle \int_0^\infty \mathds{1}_{\{Z(t) = 0\}} dt \bigg] = 0,
		$$
		for each $z \in S$, thus in order to complete the proof it suffices to prove that
		\begin{equation}\label{zero time epsilon}
			\E_\mu^z\bigg[\displaystyle \int_0^\infty \mathds{1}_{\{Z(t) \in \partial S_{i, }\}} dt \bigg] = 0, \text{ for } i =1,2.   
		\end{equation}
	 We prove this result for $i = 1$; the result then follows for  $i = 2$ by symmetry.
		
		For each $\varepsilon > 0$, define $S^{\varepsilon} \subset S$  by $S^\varepsilon = S+(\varepsilon,0)$, 
		i.e., a wedge with vertex at  $(\varepsilon, 0)$ and edges $\partial S_1^\varepsilon = \{(x,0)\in{\mathbb R}^2,\ x>\varepsilon\}$ and $\partial S_2^\varepsilon = \{(\varepsilon,0) + \lambda(\cos\xi,\sin\xi),\ \lambda> 0\}$ (recall that $\xi$ is the angle of the wedge $S$).
		
		Next we shall recursively define the $({\mathcal F}_t^z)$ stopping times $\bar\sigma_k^{\varepsilon,T},\bar\tau_k^{\varepsilon,T}$ for $k\ge 1$, for every $T>0$. We define
		
			$$
		\bar{\sigma}_1^{\varepsilon,T} = \inf\Big\{t \ge 0: Z(t) \in  S^{\varepsilon}\Big\}\wedge T,
		$$
and		
		$$
		\bar{\tau}_k^{\varepsilon,T} = \inf\Big\{t \ge \bar{\sigma}_{k}^{\varepsilon,T}: Z_t \in \partial S_2^{2\varepsilon/3}\Big\}\wedge T,\quad k \ge 1,
		$$
			$$
		\bar{\sigma}_k^{\varepsilon,T} = \inf\Big\{t \ge \bar{\tau}_{k-1}^{\varepsilon,T}: Z_t \in \partial S_2^{\varepsilon}\Big\}\wedge T,\quad k \ge 2;
		$$	
		
		Let $Z_t = (Z_1(t), Z_2(t))$. 
		Let $C>z_2$ be an  arbitrary constant. We define the $({\mathcal F}_t^z)$ stopping time 
		$$
		T_C = \inf \Big\{t \ge 0:  Z_t^{(2)} \ge C \Big \},
		$$
		and in order to simplify the notation, we also introduce the stopping times $\bar{\sigma}_k = \bar{\sigma}_k^{\varepsilon,T} \wedge T_C$ and $\bar{\tau}_k = \bar{\tau}_k^{\varepsilon,T}\wedge T_C$. Notice that for all $t\le T$, $t\in[\bar\sigma_k,\bar\tau_k]$ implies that $Z_t\in S^{2\varepsilon/3}$ and $Z_2(t)\le C$.
		
		
		Let $f_{\varepsilon,C} \in C_b^2(S)$ such that
		\begin{equation} f_{\varepsilon,C}(x,y)= 
		\begin{cases}
			0,\ {\rm if}\ (x,y)\in S\setminus S^{\varepsilon/3},\\
			y,\ {\rm if}\ (x,y)\in S^{2\varepsilon/3}, y\le C.
			\end{cases}\label{f}\end{equation}
	In addition we require that $f_{\varepsilon,C}(x,0) =0$ for all $x\ge 0$, and $D_2 f_{\varepsilon, C}\ge 0$ on $\partial S_2$. It follows from \eqref{f} that   $D_1 f_{\varepsilon, C} = 0$ on $\partial S_1$. We show in Lemma \ref{f:function} in the Appendix that such function indeed exists. By the definition of the submartingale problem
		\begin{equation}V_1=\left\{V_1(t)=f_{\varepsilon,C}(Z_t)- \int_0^t \left(\mu\cdot\nabla f_{\varepsilon,C}(Z_s)-\frac{1}{2}\triangle f_{\varepsilon,C}(Z_s)\right)ds;\ t\ge 0\right\}\label{V1:def}\end{equation}
		is a regular submartingale under $\PP_\mu^z$ on $({\mathcal F}_t^z)$, thus by Theorem 1.4.14 in \cite{karatzas2012brownian} it has a unique Doob-Meyer decomposition
		\begin{equation}V_1(t)=M(t)+A(t)\label{DM:decomposition}\end{equation}
		where $M$ is a continuous martingale and $A$ is a continuous increasing process. For the deffinition of regular submartingales see Definition 1.4.12 in \cite{karatzas2012brownian}. For an arbitrary $k\ge 1$ we have $f_{\varepsilon,C}(Z_t)=Z_2(t)$, $\mu\cdot\nabla f_{\varepsilon,C}(Z_t)=\mu_2$, and $\triangle f_{\varepsilon,C} (Z_t)=0$ whenever $t\in[\bar\sigma_k,\bar\tau_k]$, hence 
	by \eqref{DM:decomposition} and \eqref{V1:def}
	\begin{equation} Z_2(t)=Z_2(\bar\sigma_k) +M(t)-M(\bar\sigma_k) +A(t) - A(\bar\sigma_k) +\mu_2(t-\bar\sigma_k)\label{shifted}\end{equation}
	for $t\in[\bar\sigma_k,\bar\tau_k]$.
	Next we are going to establish the following two properties. The first is that
	\begin{equation} \int_{\bar \sigma_k}^{\bar\tau_k}1_{\{Z_2(t)>0\}}dA(t)=0,\quad \PP_\mu^z{\rm -a.s.},\label{first}\end{equation}
	and the second is that
	\begin{equation}
		\int_{\bar\sigma_k}^{\bar\tau_k} d\left(\langle M\rangle_t-t\right)=0,\quad, \PP_\mu^z{\rm-a.s.}
		\label{second}\end{equation}
We start with proving \eqref{first}. For any  $\delta>0$ and $k\ge 1$ we define a sequence of $({\mathcal F}^z_t)$ stopping times
$$\theta_1^\delta=\inf\{t\ge \bar\sigma_k:\ Z_z(t)\ge\delta\}\wedge\bar\tau_k, \quad \vartheta_1^\delta=\inf \left\{t\ge \theta_1^\delta:\ Z_2(t)=\frac {\delta}{2}\right\}\wedge\bar\tau_k,$$		
	$$	\theta_n^\delta=\inf \left\{t\ge \vartheta_{n-1}^\delta:\ Z_2(t)= {\delta}\right\}\wedge\bar\tau_k,\quad n\ge 2,$$
$$	\vartheta_n^\delta=\inf \left\{t\ge \theta_n^\delta:\ Z_2(t)=\frac {\delta}{2}\right\}\wedge\bar\tau_k,\quad n\ge 2.$$
Notice that $[\theta_n^\delta,\vartheta_n^\delta]$ is a sub-interval of $[\bar\sigma_k,\bar\tau_k]$ such that for $t\in  [\theta_n^\delta,\vartheta_n^\delta]$ we have $Z_2(t)\ge\delta/2$. All these stopping times are finite because by definition $\bar\tau_k\le T$. Let $g_1\in C_b^2({\mathbb R})$ be an arbitrary function such that $g_1'(0)=0$, and $g_1(x)=x$ whenever $x\ge \delta/2$. Relation $g_1'(0)=0$ implies that 
$$V_2=\left\{ V_2(t)=g_1\left(f_{\varepsilon,C}(Z_t)\right)-\int_0^t \left( \mu\cdot \nabla (g_1\circ f_{\varepsilon,C})(Z_s)- \frac{1}{2} \triangle  (g_1\circ f_{\varepsilon,C})(Z_s) \right)ds,\ t\ge 0\right\}$$
is a martingale with respect to the filtration $({\mathcal F}_t^z)$ under $\PP_\mu^z$. 
\sloppy Therefore, 
 $\left\{V_2\left((t\vee\theta_n^\delta)\wedge\vartheta_n^\delta\right),\ t\ge 0\right\}$
is also a martingale  with respect to the filtration \break $\left\{{\mathcal F}_{(t\vee\theta_n^\delta)\wedge\vartheta_n^\delta}^z,\ t\ge 0\right\}$ under $\PP_\mu^z$. For all $s\in[\theta_n^\delta,\tau_n^\delta]$ we have $g_1(f_{\varepsilon,C}(Z(s)))=Z_2(s)$,
$\mu\cdot\nabla(g_1\circ f_{\varepsilon,C})(Z(s))=\mu_2$ and $\triangle(g_1\circ f_{\varepsilon,C})(Z(s))=0$, hence 
 for all $t\ge 0$
$$V_2\left((t\vee\theta_n^\delta)\wedge\vartheta_n^\delta\right)= V_2(\theta_n^\delta)-Z_2(\theta_n^\delta) + Z_2\left((t\vee\theta_n^\delta)\wedge\vartheta_n^\delta\right) -\mu_2\left((t\vee\theta_n^\delta)\wedge\vartheta_n^\delta-\theta_n^\delta\right)$$
thus $\left\{Z_2\left((t\vee\theta_n^\delta)\wedge\vartheta_n^\delta\right) -\mu_2\left((t\vee\theta_n^\delta)\wedge\vartheta_n^\delta\right),\ t\ge 0\right\}$	is also  a martingale with respect to the filtration $\left\{{\mathcal F}^z_{(t\vee\theta_n^\delta)\wedge\vartheta_n^\delta},\ t\ge 0\right\}$ under $\PP_\mu^z$. On the other hand, from \eqref{shifted} follows that 
$$	Z_2\left((t\vee\theta_n^\delta)\wedge\vartheta_n^\delta\right) -\mu_2\left((t\vee\theta_n^\delta)\wedge\vartheta_n^\delta\right)=$$
$$Z_2(\theta_n^\delta) -\mu_2\theta_n^\delta +M( (t\vee\theta_n^\delta)\wedge\vartheta_n^\delta  )-M(\theta_n^\delta) +A( (t\vee\theta_n^\delta)\wedge\vartheta_n^\delta ) - A(\theta_n^\delta) ,$$
for all $t\ge 0$.
However, the left-hand side  in the above identity is a martingale with respect to the filtration $\left\{{\mathcal F}^z_{(t\vee\theta_n^\delta)\wedge\vartheta_n^\delta},\ t\ge 0\right\}$ under $\PP_\mu^z$, and so is $M( (t\vee\theta_n^\delta)\wedge\vartheta_n^\delta  )$ on the right hand side ($t\ge 0$). Therefore, $A$ must be constant on $[\theta_n^\delta,\tau_n^\delta]$, $\PP_\mu^z$-a.s. This holds for all $n\ge 1$, hence
\begin{equation}\int_{\hat\sigma_k}^{\hat\tau_k} \sum_{n=1}^\infty 1_{[\theta_n^\delta,\vartheta_n^\delta]}(t)dA(t)=0,\quad \PP^z_\mu{\rm -a.s.}\label{almost:there}\end{equation}
If $t\in[\bar\sigma_k,\bar\tau_k]$ and $Z_2(t)>\delta$, then $Z_2(t)\in[\theta_n^\delta,\vartheta_n^\delta]$ for some $n\ge 1$, hence by \eqref{almost:there}
$$\int_{\hat\sigma_k}^{\hat\tau_k} 1_{\{Z_2(t)>\delta\}}(t)dA(t)=0,\quad \PP_\mu^z{\rm-a.s.,}$$
and \eqref{first} follows.

Next we are going to show \eqref{second}. Let $g_2\in C_b^2({\mathbb R})$ arbitrary such that $g_2(x)=x^2$ whenever $|x|\le C$. Then
\begin{equation}V_3=\left\{V_3(t)= g_2\left(f_{\varepsilon,C}(Z_t)\right)-\int_0^t \left( \mu\cdot \nabla (g_2\circ f_{\varepsilon,C})(Z_s)- \frac{1}{2} \triangle  (g_2\circ f_{\varepsilon,C})(Z_s) \right)ds,\ t\ge 0\right\}\label{U:def}\end{equation}
is a martingale under $\PP^z_\mu$ with respect to the filtration $({\mathcal F}^z_t)$, and for $t\in[\bar\sigma_k,\bar\tau_k]$ we have  $g_2\left(f_{\varepsilon,C}(Z_t)\right)=\left(Z_2(t)\right)^2$, $\mu\cdot \nabla (g_2\circ f_{\varepsilon,C})(Z_t)=2\mu_2 Z_2(t)$ and $\triangle  (g_2\circ f_{\varepsilon,C})(Z_t)=2$, hence by Ito's rule applied to $g_2(f_{\varepsilon,C}(Z_t))$
and by \eqref{shifted}
$$Z_2^2(t)=Z_2^2(\bar\sigma_k) +\int_{\bar\sigma_k}^t2Z_2(s)dM(s)  + \int_{\bar\sigma_k}^t 2\mu_2Z_2(s)ds +\langle M\rangle_t - \langle M\rangle_{\bar\sigma_k},$$
for $t\in[\bar\sigma_k,\bar\tau_k]$. We note that the $\int_{\bar\sigma_k}^t2Z_2(s) dA(s)$ term vanished because of \eqref{first}.
From this and from \eqref{U:def} follows that
$$V_3(t)= V_3(\bar\sigma_k) +\int_{\bar\sigma_k}^t2Z_2(s)dM(s)  + \int_{\bar\sigma_k}^t d\left(\langle M\rangle_s -s\right),$$
for $t\in[\bar\sigma_k,\bar\tau_k]$. The process $\left\{V_3\left((t\vee\bar\sigma_k)\wedge\bar\tau_k\right),\ t\ge 0\right\}$ is a martingale with respect to  the filtration $\left\{{\mathcal F}_{(t\vee\theta_n^\delta)\wedge\vartheta_n^\delta},\ t\ge 0\right\}$ under $\PP_\mu^z$, and can be  written by substituting $(t\vee\bar\sigma_k)\wedge\bar\tau_k$ for $t$ in the above identity as
$$V_3(t\vee\bar\sigma_k)\wedge\bar\tau_k)= V_3(\bar\sigma_k) +\int_{\bar\sigma_k}^{t\vee\bar\sigma_k)\wedge\bar\tau_k  }
2Z_2(s)dM(s) + \int_{\bar\sigma_k}^{  t\vee\bar\sigma_k)\wedge\bar\tau_k }   d\left(\langle M\rangle_s -s\right),$$
for all $t\ge 0$.
Since the left-hand side is a martingale  with respect to  the filtration $\left\{{\mathcal F}^z_{(t\vee\theta_n^\delta)\wedge\vartheta_n^\delta},\ t\ge 0\right\}$ under $\PP_\mu^z$, \eqref{second} follows. 
The by \eqref{shifted}, $Z_2$ is a 1-dimensional Brownian motion with drift $\mu_2$ reflected at zero in $[\bar\sigma_k,\bar\tau_k]$.
Therefore,
$$\int_{\bar\sigma_k}^{\bar\tau_k}1_{\{Z_2(t)=0\}}dt=0,\quad \PP^z_\mu\text {-a.s.}$$
This holds for every $k\ge 1$, hence we also have
$$\sum_{k=1}^\infty\int_{\bar\sigma_k}^{\bar\tau_k}1_{\{Z_2(t)=0\}}dt=0,\quad \PP^z_\mu\text {-a.s.,}$$
and from this
$$\int_0^{T\wedge T_C} 1_{\{Z_1(t)\ge \varepsilon, Z_2(t)=0\}}dt=0, \quad \PP^z_\mu\text {-a.s.}$$
The last identity follows because $t\le T\wedge T_c$, $Z_1(t)\ge\varepsilon$ and $Z_2(t)=0$ implies that $t\in [\bar\sigma_k,\bar\tau_k]$ for some $k\ge 1$. The statement of the Lemma now follows by $T,C\uparrow\infty$ and $\varepsilon\downarrow 0$.

	\end{proof}
	
	Let $\{\PP_\mu^z,z\in S\}$ be a solution to the submartingale problem with a drift $\mu$.  Next we shall create a process $X$  which is a Brownian motion with drift $\mu$ started at $z$ on $(C_S,{\mathcal M},({\mathcal M}_t),\PP_\mu^z)$, for every $z\in S$. We already know that such process exists for the solution that we created in Proposition \ref{existence theorem}. However, for proving the uniqueness of the solution, we need to show the existence of such process $X$ for every solution of the submartingale problem.
	Such construction has been carried out in \cite{lakner2019roughness} and in \cite{kang2017submartingale} for the case of zero drift. The generalization to the case of non-zero drift requires only a few obvious changes to the proofs in the case of zero drift, so here we shall only state the results (Lemmas \ref{w:delta} and \ref{construct the Brownian motion using the convergence})  without proofs.
	
	For each $\delta > 0$, let $S_\delta \subset S$ be the closed set defined in 
	the complex plane by $S_\delta = S+\delta e^{i\xi/2}$. So $S_\delta$ is a wedge with vertex at $\delta\left(\cos\left(\xi/2\right),\sin\left(\xi/2\right)\right)$, such that it is included in $S$ and has edges parallel with the respective edges of $S$. 
	
	 Set $\tau_0^\delta = 0$, and, for each $k \ge 1$, recursively define
	$$
	\sigma_k^\delta = \inf\{ t \ge \tau_{k-1}^\delta: Z(t) \in S_{2\delta}\} \text{ and } \tau_k^\delta = \inf\{t \ge \sigma_k^\delta: Z(t) \in S\backslash S_{\delta}\}.
	$$
	By Problem 1.2.7  in Karatzas and Shreve \cite{karatzas2012brownian}, $\sigma_k^{\delta}$ and $\tau_k^{\delta}$  are stopping times relative to {$\{\calM_t, t\ge 0\}$} for every $k\ge 1$. 
	
	For each $k \ge 1$ and $\delta > 0$, define the process $\{W_{(k)}^\delta(t), t \ge 0\}$ by setting
	$$
	W_{(k)}^\delta(t) = Z(t \wedge \tau_k^\delta) - Z(t \wedge \sigma_k^\delta) - (\ttkd - \tskd)\mu, \quad t\ge 0, 
	$$
	and then define the process $\{W^\delta(t),t\ge 0\}$ by setting 
	$$
	W^\delta(t) = \sum_{k=1}^{\infty} W_{(k)}^\delta(t),\quad t\ge 0. 
	$$

\begin{lem}\label{w:delta}  For every $\delta>0$ and $z\in S$ the process $W^\delta$ is a square-integrable martingale on $(C_S, {\mathcal M}, ({\mathcal M}_t), \PP^z_\mu)$.
	\end{lem}

	\begin{lem}\label{construct the Brownian motion using the convergence}
		 There exists a process $W$ on {$(C_S, \calM, \calM_t)$} such that for every $z\in S$  it is a standard $2$-dimensional Brownian motion under $\PP_\mu^z$ starting at zero, and 
		for every fixed $T>0$ we have
		\begin{equation}
			E_\mu^z\left[\|W(T)-W^\delta(T)\|^2\right]\to 0,\quad{\rm as}\ \delta\to 0.\label{goestozero}
		\end{equation}
		\end{lem}

	Next we shall define the process $X$ by
		\begin{equation}
		X(t, \omega) = \omega(0) + W(t, \omega)+\mu t,\quad t\ge 0.\label{Xdef}
	\end{equation}
 We  define the process $Y$ by
 \begin{equation}Y(t) = Z(t)-X(t), \quad  t\ge 0.\label{Y:def}\end{equation}
	We shall say that a function $f:{\mathbb R}_+\mapsto {\mathbb R}^2$ is {\it flat} on an interval $[s,t]\subset {\mathbb R}_+$, if for every $u\in [s,t]$ we have $f(u)=f(s)$. 
	
	\begin{lem}\label{Existence of BM under drift measure}
		Let $\{\PP^z_\mu; z\in S\}$ be an arbitrary solution of the submartingale problem, and
		let $X$ and $Y$ be the processes defined above. Then the following two statements hold for every $z\in S$:
		\begin{itemize}
	\item[1.]	 Under $\PP_\mu^z$ the process $X$ is a standard 2-dimensional Brownian motion on {$(C_S, \calM, \calM_t)$} with drift $\mu$, started at $z$;
\item[2.]		 for every $n\in{\mathbb N}_+$, and $\delta>0$, the sample paths of $Y$ are flat on $[\sigma_n^\delta, \tau_n^\delta]$, $\PP_\mu^z$-a.s.
		 \end{itemize}
	\end{lem}
	\begin{proof} 
The first statement follows 	from Lemma \ref{construct the Brownian motion using the convergence} and property 1 in the definition of the Submartingale Problem. Next we shall prove the second statement.
By the definition of $w^\delta$, the sample paths of 
\begin{equation}\{Z_t-w^\delta(t)-\mu t,\ t\ge 0\}\ \text{are flat on}\ [\sigma_n^\delta,\tau_n^\delta],\label{flat}\end{equation}
for each $\delta>0$, $n\ge 1$. On the other hand, for every $\delta>0$, $n\ge 1$ there exists $k\ge 1$ (depending on the sample path) such that $[\sigma_n^\delta,\tau_n^\delta]\subset [\sigma_k^{\delta /2},\tau_k^{\delta/2}]$. This implies that the sample paths of $\{Z_t-w^{\delta/2}(t)-\mu t,\ t\ge 0\}$ are also flat on $[\sigma_n^\delta,\tau_n^\delta]$. Iterating this we get that for every $m\ge 1$ the sample paths of  $\{Z_t-w^{\delta/2^m}(t)-\mu t,\ t\ge 0\}$ are also flat on $[\sigma_n^\delta,\tau_n^\delta]$. Comparing this with \eqref{flat} we conclude that the sample paths of  $w^\delta-w^{\delta/2^m}$ are also flat on $[\sigma_n^\delta,\tau_n^\delta]$, thus for every $t\ge 0$
$$\int_0^t 1_{[\sigma_n^\delta,\tau_n^\delta]}(s) dw^{\delta/2^m}(s)= \int_0^t 1_{[\sigma_n^\delta,\tau_n^\delta]}(s) dw^{\delta}(s).$$
Taking limit as $m\to \infty$ and using  \eqref{goestozero} we get that 
$$\int_0^t 1_{[\sigma_n^\delta,\tau_n^\delta]}(s) dw(s) = \int_0^t 1_{[\sigma_n^\delta,\tau_n^\delta]}(s) dw^{\delta}(s),$$
$\PP_\mu^z$-a.s.
This identity and \eqref{flat} imply that 	the sample paths of $\{Z_t-w(t)-\mu t,\ t\ge 0\}$ are flat on  $[\sigma_n^\delta,\tau_n^\delta]\cap[0,t]$, $\PP^z_\mu$-a.s. Since $t\ge 0$ was arbitrary,  this and \eqref{Xdef} imply what we wanted to prove.

	\end{proof}
	
	\begin{lem}\label{unique:main}
	Suppose that $Q_1$ and $Q_2$ are mutually absolutely continuous probability measures on ${\mathcal M}$, both satisfying properties 1,2, and 3 of Definition \ref{defn of the submartingale problem} with $\PP^z_\mu$ replaced by $Q_i$ ($i=1,2$). Then there exist probability measures $\tilde Q_i$ on ${\mathcal M}$ for $i=1,2$, such that conditions 1,2, and 3 of Definition  \ref{defn of the submartingale problem} are satisfied with $P_\mu^z$ replaced by $\tilde Q_i$  and $\mu$ replaced by the zero vector.
	Furthermore,  for every $T\ge 0$ there exist probability measures $\tilde Q^T_1$ and $\tilde Q^T_2$ on ${\mathcal M}$ such that for all $T\ge 0$, $A\in{\mathcal M}_T$, and $i=1,2$ we have
	  $\tilde Q_i^T(A)=\tilde Q_i(A)$, $\tilde Q_i^T$ and $Q_i$ are mutually absolutely continuous, and
	 \begin{equation}
	 	\frac{d\tilde Q_1^T}{dQ_1} = 	\frac{d\tilde Q_2^T}{dQ_2}.\label{abscont}
	\end{equation}
	
	\end{lem}
	
	\begin{proof}
		Let $Q_1$ and $Q_2$ as above, and
		let $(X^i, Y^i)$ be as in Lemma \ref{Existence of BM under drift measure} defined under $\Q_i$. Since $X^i$ is defined by $L^2(\Q_i)$ convergence, this implies that $(X^1,Y^1) = (X^2, Y^2)$, which we shall from here on denote by $(X, Y)$.
In this proof we shall use the Doob-Meyer decomposition which requires that the probability space satisfies the ``usual conditions", and for this purpose we have to augment the probability space $(C_S,{\mathcal M},({\mathcal M}_t), Q_i)$, $i=1,2$; let this augmentation be $(C_S, {\mathcal F},({\mathcal F}_t),Q_i)$.  The measures $Q_1$ and $Q_2$ are mutually absolutely continuous, hence the filtration $({\mathcal F}_t)$ and the sigma field ${\mathcal F}$ do not depend on $i=1,2$.
 For technical details concerning the augmentation of a probability space please see Remark \ref{aug} in the Appendix. In this proof all processes live on the augmented space $(C_S, {\mathcal F}, ({\mathcal F}_t))$, unless specified otherwise. 
 
 Next, note that for each $T>0$ and $\delta > 0$ and $n \geq 1$, 
		$$
		\Big\{Z(\tau_n^\delta \wedge t) - Z(\sigma_n^\delta \wedge t),\  t \in [0, T]\Big\}
		$$ 
		is a semimartingale under both $Q_1$ and $Q_2$ with respect to the filtration $({\mathcal F}_t)$.  Indeed, it can be written 
by Lemma \ref{Existence of BM under drift measure}  and by \eqref{Y:def}  as
		$$Z(\tau_n^\delta \wedge t) - Z(\sigma_n^\delta \wedge t)= X(\tau_n^\delta \wedge t) - X(\sigma_n^\delta \wedge t)
		, \quad t \in [0,T] .
		$$
		
		Now let $f \in C_b^2(S)$ such that $D_if\ge 0$ on $\partial S_i$ for $i=1,2$, and $f$ is constant in a neighborhood of the origin. Then, by It\^{o}'s rule we have that for $t \in [0,T]$,
		\begin{equation}
			f(Z(  t\wedge\tau_n^\delta  ) )  = f(Z(t\wedge \sigma_n^\delta)) +  
			 \int_{t\wedge \sigma_n^\delta}^{t\wedge \tau_n^\delta} \nabla f(Z(s)) dX(s) +  
			 \frac{1}{2} \int_{t\wedge \sigma_n^\delta}^{t\wedge\tau_n^\delta} \Delta f(Z(s)) ds, \label{decomposition1}
\end{equation}
	$Q_i$-a.s., $i=1,2$.	On the other hand, by condition 2 of Definition \ref{defn of the submartingale problem} and by Theorems I.4.10 and I.4.14 in   \cite{karatzas2012brownian}, we have for $i=1,2$, the unique Doob-Meyer decomposition
		\begin{equation}
			f(Z(t)) = f(z) + \int_0^t \nabla f(Z(s))\cdot \mu ds + \frac{1}{2} \int_0^t \Delta f(Z(s)) ds + M^i(t) + A^i(t),\ t\le T, \label{doob::meyer:f}
		\end{equation}
		where $M^i$ is a continuous  martingale and $A^i$ is a continuous, increasing process on $(C_S,{\mathcal F}, ({\mathcal F}_t), Q_i)$, with $M^i(0) = A^i(0) = 0$. By Proposition 16.32 in Bass \cite{bass2011stochastic}, we also have that for $T \geq 0$,
		$$
		\E^{\Q_i} [||M^i(T)||^2] < \infty,~i=1,2.
		$$
		
		Let $W$ as in \eqref{Xdef}, that is, $W(t)=X(t)-z-\mu t$,
		and for $i=1,2,$ let $S^i(W)$ be the class of $\R^2$-valued processes on $(C_S,{\mathcal F}, ({\mathcal F}_t))$ such that $U \in S^i(X)$ if it has the form 
		$$
		U(t) = \int_0^t G(s)dW(s),~t \in [0,T],
		$$
		for some $2$-dimensional process $G$ such that 
		
		$$
		\E^{\Q_i} \bigg[\int_0^T ||G(s)||^2 ds \bigg] < \infty.
		$$
		Then, by Theorem IV.36 and Corollary 1 to Theorem IV.37 in \cite{protter2005stochastic}, there exists a $\R^2$-valued process $H^i$ such that 
		$$
		M^i(t) = \int_0^t H^i(s)dW(s) + N^i(t),~t \in [0,T],
		$$
		where
		\begin{itemize}
			\item [(i)] $N^i$ is a square-integrable martingale under $\Q_i$, 
			\item [(ii)] $N^i$ is strongly orthogonal to every member of $S^i(W)$ under $\Q_i$, that is, $N^iU$ is a $\Q_i$-martingale for each $U \in S^i(W)$, 
			\item [(iii)] $\E^{\Q_i} \bigg[ \displaystyle \int_0^T ||H^i(s)||^2 ds \bigg] < \infty$.
		\end{itemize}
		Now, by \eqref{doob::meyer:f}  we have for $t \in [0,T]$,
		\begin{eqnarray}
			f(Z(t)) &=& f(z) + \int_0^t \nabla f(Z(s)) \cdot \mu ds + \frac{1}{2} \int_0^t \Delta f(Z(s)) ds \nonumber \\
			&+& \int_0^t H^i(s)d W(s) + N^i(t) + A^i(t),\ Q_i-{\rm a.s.},\label{intermediate}
		\end{eqnarray}
		
		hence, by \eqref{intermediate}  we have for $t \in [0,T]$,
		$$
			f(Z(t\wedge\tau_n^\delta))  = f(Z(t\wedge\sigma_n^\delta)) 
			+ \int_{t\wedge \sigma_n^\delta}^{t\wedge \tau_n^\delta} \nabla f(Z(s))\cdot \mu ds + {1\over 2}\int_{t\wedge \sigma_n^\delta}^{t\wedge\tau_n^\delta} \Delta f(Z(s)) ds$$
			\begin{equation}+ \int_{t\wedge \sigma_n^\delta}^{t\wedge\tau_n^\delta} H^i(s)d W(s) + N^i(t\wedge \tau_n^\delta) - N^i(t\wedge\sigma_n^\delta) 
			+ A^i(t\wedge \tau_n^\delta) - A^i(t\wedge\sigma_n^\delta), \label{decomposition2}
		\end{equation}
		$Q_i$-a.s. Now, for each $i=1,2$, we have two Doob-Meyer decompositions of the submartingale 
		$$\{f(Z(t\wedge\tau_n^\delta))  - f(Z(t\wedge\sigma_n^\delta))-
		\int_{t\wedge \sigma_n^\delta}^{t\wedge \tau_n^\delta} \nabla f(Z(s))\cdot \mu ds - {1\over 2}\int_{t\wedge \sigma_n^\delta}^{t\wedge\tau_n^\delta} \Delta f(Z(s)) ds
		,\ t\in[0,T]\},$$ 
		 \eqref{decomposition1} and \eqref{decomposition2}.  Hence, by the uniqueness of the Doob-Meyer decomposition, for each $t \in [0,T]$
		\begin{eqnarray}
			N^i(t\wedge \tau_n^\delta) - N^i(t\wedge\sigma_n^\delta) + \int_{t\wedge \sigma_n^\delta}^{t\wedge\tau_n^\delta} H^i(s)d W(s) = \int_{t\wedge \sigma_n^\delta}^{t\wedge \tau_n^\delta} \nabla f(Z(s)) dW(s), \label{d-m}
		\end{eqnarray}
		 $Q_i$- a.s. $N_i$ is strongly orthogonal to every member of $S^i(W)$, and from \cite{protter2005stochastic}, Theorem IV.37  follows that $ 	\{N^i(t\wedge \tau_n^\delta) - N^i(t\wedge\sigma_n^\delta),\ t\in[0,T]\}  $ is also strongly orthogonal to every member of $S^i(W)$ under $Q_i$.
		 However, by the above relation it is also a member of $S^i(W)$, hence it follows that
		$N^i(t\wedge \tau_n^\delta) - N^i(t\wedge\sigma_n^\delta) = 0$ for $t\in[0,T]$. Then by \eqref{d-m} we also have
		$$
		\E^{\Q_i} \bigg[ \displaystyle \int_{t\wedge \sigma_n^\delta}^{t\wedge \tau_n^\delta} ||H^i(s) - \nabla f(Z(s))||^2 ds \bigg] = 0.
		$$
		Hence, 
		\begin{eqnarray}\label{ineq1}
			&&\E^{\Q_i} \bigg[ \int_0^t ||H^i(s) - \nabla f(Z(s))||^2ds\bigg] \nonumber \\
			&=& \E^{\Q_i} \bigg[ \int_0^t \sum_{n=2}^{\infty} \mathds{1}_{\big\{s \in [\tau_{n-1}^\delta, \sigma_n^\delta]\big\} }||H^i(s) - \nabla f(Z(s))||^2ds\bigg] \nonumber \\
			&\le& \E^{\Q_i} \bigg[ \int_0^t \mathds{1}_{\big\{Z(s) \in S_{2\delta}^c\big\}}||H^i(s) - \nabla f(Z(s))||^2ds\bigg],
		\end{eqnarray}
		where $S_{2\delta}^c = S\backslash S_{2\delta}$. Moreover, by the dominated convergence theorem, 
		$$
		(\ref{ineq1}) \ra \E^{\Q_i} \bigg[ \int_0^t \mathds{1}_{\big\{Z(s) \in \partial S\big\}}||H^i(s) - \nabla f(Z(s))||^2 ds\bigg] = 0~\mathrm{as}~\delta \rightarrow 0,
		$$
		where the last identity is by Lemma \ref{Occupation:time} . By \eqref{intermediate}, now follows that for $t \in [0,T]$,
		\begin{eqnarray}
			f(Z(t)) &=& f(z) + \int_0^t \nabla f(Z(s)) \cdot \mu ds + \frac{1}{2} \int_0^t \Delta f(Z(s)) ds \nonumber \\
			&+& \int_0^t \nabla f(Z(s)) dW(s) + N^i(t) + A^i(t),\ Q_i{\rm -a.s.}\label{almost}
		\end{eqnarray}
		
		Next, for each $t \geq 0$ let
		$$
		\tilde\zeta(t) = \exp\bigg\{-\mu \cdot (X(t)-z) + \frac{1}{2}||\mu||^2t\bigg\},
		$$
		and for each $i=1,2$, and $T \geq 0$ define the measure $\tilde{\Q}_i^T$ by setting
		\begin{eqnarray}
			\frac{d \tilde{\Q}_i^T}{d \Q_i} &=& \tilde\zeta(T). \label{def:rn:deriv}
		\end{eqnarray}
		Then, under $\tilde{\Q}_i^T$, $\{ X(t), t\in [0,T]\} $ is a Brownian motion (without drift) started at $z$, and by \eqref{almost} we have 
		\begin{eqnarray}
			f(Z(t)) &=& f(z) + \int_0^t \nabla f(Z(s)) dX(s)  \nonumber \\
			&+& \frac{1}{2} \int_0^t \Delta f(Z(s)) ds +  N^i(t) + A^i(t),\ Q_i{\rm -a.s.}\label{decomposition}
		\end{eqnarray}
		However, note that $N^i$ is also a $\tilde{\Q}_i^T$-martingale on $[0,T]$. 
		Indeed, by \eqref{def:rn:deriv}, $N^i$ is a $\tilde{\Q}_i^T$-martingale if $N^i\tilde\zeta$ is a $\Q_i$-martingale on $[0,T]$. But this follows since from the fact that $N_i$ is strongly orthogonal to every member of $S^i(X)$ under $\Q_i$, and  by its definition $\tilde\zeta-1 \in S^i(W)$ under $\Q_i$. Just like in the proof of Proposition \ref{existence theorem}, there exists a probability measure $\tilde \Q_i$ on ${\mathcal M}$ such that $\tilde \Q_i(A)=\tilde \Q_i^T(A)$, whenever $A\in{\cal M}_T$. Thus, by \eqref{decomposition} both $\tilde \Q_1$ and $\tilde \Q_2$ satisfy property 2 of Definition \ref{defn of the submartingale problem} with $\PP^z_\mu$ replaced by either $\tilde Q_1$ or $\tilde Q_2$, and $\mu$ replaced by the zero vector. 
$\tilde Q_i^T$ and $Q_i$ are mutually absolutely continuous because of \eqref{def:rn:deriv} and because $\tilde\zeta(T)>0$, a.s. under $Q_i$. Relation \eqref{abscont} also follows from \eqref{def:rn:deriv}. Properties 1 and 3 of  Definition \ref{defn of the submartingale problem} are satisfied if $\PP^z_\mu$ is replaced by $\tilde Q_i$ because they are satisfied if we replace $\PP^z_\mu$ by $Q_i$, and we already established that $\tilde Q_i^T$ and $Q_i$ are mutually absolutely continuous. Property 1 follows immediately from this. Property 3 can be shown by first showing it with the $\infty$ in the upper limit of the integral replaced by $T$, then taking $T\to\infty$.

	\end{proof}
	
\noindent	{\it Proof of Theorem \ref{Exsistence and Uniqueness of Submartingale Problem} .} In light of  Propositions \ref{existence theorem}  the only missing part is the proof of uniqueness.
	Let $z \in S$ and suppose that $\PP_1^z$ and $\PP_2^z$ are two probability measures satisfying conditions 1, 2 and 3 in Definition \ref{defn of the submartingale problem} of the submartingale problem with drift.  Let 
	$$
	\Q_1 = \frac{1}{3} \PP_1^z + \frac{2}{3} \PP_2^z \text{\quad and \quad} \Q_2 = \frac{2}{3} \PP_1^z+ \frac{1}{3} \PP_2^z.
	$$
	Then, one can check that each $\Q_i$ also satisfies conditions 1,2, and 3 in 
	Definition \ref{defn of the submartingale problem} of the submartingale problem with drift. In addition, $\Q_1$ and $\Q_2$ are mutually absolutely continuous. In order to complete the proof, it is therefore sufficient to show that $\Q_1 \equiv \Q_2$. 
	By Lemma \ref{unique:main} there exist probability measures $\tilde Q_i$, $i=1,2$, such that properties 1,2, and 3 of Definition \ref{defn of the submartingale problem} are satisfied with $P^z_\mu$ replaced by $\tilde Q_i$, and $\mu$ replaced by the zero vector. 
	  The uniqueness result in Section 3.1 of \cite{varadhan1985brownian} implies  $\tilde \Q_1=\tilde \Q_2$. Using the probability measures $\tilde Q_i^T$ from the same proposition, we have now that  
	   $$\tilde Q_1^T|_{{\mathcal M}_T} =    \tilde Q_2^T |_{{\mathcal M}_T} .$$ From \eqref{abscont}
 follows that 
	    $$ Q_1|_{{\mathcal M}_T} =     Q_2 |_{{\mathcal M}_T} .$$
	    Since $T$ was arbitrary, $Q_1=Q_2$ follows.\hfill$\square$
	 
	 \vskip .2in

\begin{proof}[Proof of Theorem \ref{variation}] Recall from the proof of Proposition \ref{existence theorem} that for every $T\ge 0$ there exists the probability measure $\PP^z_{\mu,T}$ which is mutually absolutely continuous with respect to $\PP^z_0$, and coincides with $\PP^z_\mu$ on ${\mathcal M}_T$. By Theorem 2.6 in \cite{lakner2019roughness}, formulas \eqref{var_finite} and \eqref{var_infinite} hold for $\mu=0$. Then by the mutual absolute continuity of $\PP_0^z$ and $\PP^z_{\mu,T}$, \eqref{var_finite} and \eqref{var_infinite} also hold with $\PP_\mu^z$ replaced by $\PP_{\mu,T}^z$. Since $\PP_\mu^z$ and $\PP_{\mu,T}^z$ coincide on ${\mathcal M}_T$, both formulas follow.
    
\begin{proof}[Proof of Theorem \ref{dirichlet}]
This follows from Theorem 2.4 in \cite{lakner2019roughness} using the measure $\PP^z_{\mu,T}$ for every $T>0$, just like in the proof of Theorem \ref{variation}.
\end{proof}

\end{proof}	 
\begin{proof}[Proof of Theorem \ref{notsemimartingale}]
Suppose that $1\le\alpha<2$, and $Z$ is a continuous semimartingale on $(C_S,{\mathcal M}, {\mathcal M}_t, \PP^z_\mu)$ for some $z\in S$. Then there exists a decomposition
\begin{equation} Z(t) = z+M(t) + A(t),\quad t\in[0,\infty),\label{z:decomp}\end{equation}
where $M$ is a continuous local martingale and $A$ is a finite variation (FV) process on $(C_S,{\mathcal M}, {\mathcal M}_t, \PP^z_\mu)$ (see \cite{protter2005stochastic}, Corollary to Theorem II.31). We know from the proof of Proposition \ref{existence theorem} that for every $T\in[0,\infty)$ there exists a probability measure $P^z_{\mu,T}$ on ${\mathcal M}$ which is mutually absolutely continuous with respect to $\PP_0^z$, and $\PP^z_{\mu,T}(A) = \PP_\mu^z(A)$ for all $A\in{\mathcal M}_T$. In addition,
$$\frac{d\PP_0^z}{d\PP^z_{\mu,T}}=\frac{1}{\zeta(T)}$$
where $\zeta$ is defined under \eqref{rn_derivative}. 
We cast \eqref{z:decomp} in the form 
\begin{equation} Z(t) = z+\tilde M(t) +\tilde A(t),\quad t\in[0,\infty),\label{z:newdecomp}\end{equation}
where 
$$\tilde M(t)= M(t) -\int_0^t\zeta(s) d\left[\frac{1}{\zeta},M\right]_s$$
$$\tilde A(t)= A(t) +\int_0^t\zeta(s) d\left[\frac{1}{\zeta},M\right]_s.$$
By the Girsanov-Meyer theorem ( \cite{protter2005stochastic}, Theorem III.31) $\tilde M$ is a local martingale on $[0,T]$ and $\tilde A$ is a FV process under $\PP^z_0$. But this implies that $\tilde M$ is a local martingale on $[0,\infty)$ under $\PP_0^z$, hence $Z$ must be a semimartingale under $\PP_0^z$, which is in contradiction with the result of \cite{williams1985reflected}, Theorem 5.

\end{proof}
	 
\noindent	 {\it Proof of Theorem \ref{nonexistence}.} Suppose that $\alpha\ge 2$,  let $z \in S$ and suppose that $\PP^z_\mu$ is a probability measure on ${\mathcal M}$  satisfying properties 1, 2 and 3 in Definition \ref{defn of the submartingale problem} of the submartingale problem with drift. Selecting $Q_1=Q_2=\PP^z_\mu$ in Lemma \ref{unique:main}, it follows that there exists a probability measure $\tilde Q$ on ${\mathcal M}$ such that properties 1,2, and 3 of Definition \ref{defn of the submartingale problem} are satisfied with $\PP^z_\mu$ replaced by $\tilde Q$, and $\mu$ replaced by the zero vector. However, this is in direct contradiction with Theorem 3.11 in \cite{varadhan1985brownian}.\hfill$\square$

	\section{Proof of Theorems \ref{VW's Feller} and \ref {Strong Markov for Submartingale Problem}  }\label{section_Feller}
	First we shall prove Theorem \ref{VW's Feller}.
	
	\begin{prop}\label{weakly realatively compactness}
		The family of probability measures $\{\PP_\mu^{z_n}\}$ is tight for any sequence $\{z_n, n\ge 1\}$ in $S$ which converges to some $z \in S$. 
	\end{prop}
	\begin{proof}
		By Theorem 2.4.10 in \cite{karatzas2012brownian}, it is sufficient to show that    
		\begin{equation}\label{limisup condition for tightness}
			\lim_{\delta \downarrow 0} \sup_{n} \PP_\mu^{z_n}(\omega: m^T(\omega, \delta) \ge \varepsilon) = 0, \text{ for any } T > 0, \varepsilon > 0.    
		\end{equation}
		
		In the above,
		$$
		m^T(\omega,\delta) = \sup_{\substack{{|t-s| \le \delta} \\ {0\le s, t \le T}}} |\omega(s) - \omega(t)|.
		$$
		Using \eqref{rn_derivative} and the Cauchy-Schwarz inequality,   
		\begin{eqnarray}\label{tightness estimate}
			&&\PP_\mu^{z_n}(\omega: m^T(\omega, \delta) \ge \varepsilon) \nonumber\\
			&=& \E^{z_n}_0 [\mathds{1}_{\{m^T(\omega,\delta) > \varepsilon\}} \exp\{\mu\cdot(X(T)-X(0)) - {1 \over 2} ||\mu||^2T\}]\nonumber\\
			&\le&  \left(\E^{z_n}_0 [ \exp\{2\mu\cdot(X(T)-X(0)) - ||\mu||^2T\}]{\PP_0^{z_n}}(m^T(\omega,\delta) > \varepsilon)\right)^{1/2}\nonumber\\
			&=& \exp\left\{\frac{1}{2}||\mu||^{2}T\right\} \left({\PP_0^{z_n}}(m^T(\omega,\delta) > \varepsilon)\right)^{1/2}.
		\end{eqnarray}
		By Theorem 3.13 in \cite{varadhan1985brownian}, $\{\PP_0^{z_n}\}$ is tight hence
		$$
		\lim_{\delta \downarrow 0} \sup_{n} \PP_0^{z_n}(\omega: m^T(\omega, \delta) \ge \varepsilon) = 0, \text{ for any } T > 0, \varepsilon > 0, 
		$$
		combining with inequality (\ref{tightness estimate}), we have (\ref{limisup condition for tightness}).
	\end{proof}

	
	\begin{proof}[Proof of Theorem \ref{VW's Feller}]
		Given Proposition \ref{weakly realatively compactness}, it only remains to show that any weak limit point $\PP_\mu^*$ of the family $\{\PP_\mu^{z_n}\}$ is a solution to the submartingale problem starting from $z$, then by the uniqueness part of Theorem
		\ref{Exsistence and Uniqueness of Submartingale Problem}, $\PP_\mu^{z_n} \Rightarrow \PP_\mu^z$ as $n \ra \infty$.
		
		It is straightforward that $\PP_\mu^*$ satisfies condition $1$ in Definition \ref{defn of the submartingale problem} (the submartingale problem), since for any $k \ge 1$ and the closed set $C_k = \{\omega \in C_S: |\omega(0)-z| \le {1 \over k}\}$, $1 = \lim \sup_n \PP_\mu^{z_n}(C_k) \le \PP_\mu^*(C_k)$ hence $\PP_\mu^*$ concentrates on $\{\omega \in C_S: \omega(0)=z\}$. The condition $2$ is also satisfied, since the submartingale property is preserved under the weak convergence. Now we prove $\PP_\mu^*$ satisfies condition $3$, we need to show that if $(z_n,n\ge 1)\subset S$, $z\in S$ such that $\lim_{n\to\infty}z_n=z$ and $\PP^{z_n}_\mu \Rightarrow \PP^*_\mu$, then
		\begin{equation} \E^*_\mu\left[\int_0^\infty \mathds{1}_{\{\left(Z(t)\right) = 0\}}dt\right]=0.\label{goal}
		\end{equation}
	Let $\epsilon>0$ and $t\ge 0$ be arbitrary. Under the local uniform topology  the event $\{w\in C_S:\ |Z_t(\omega)|<\epsilon\}$ is an open set, and  $\{w\in C_S:\ |Z_t(\omega)|\le \epsilon\}$ is a closed set. By the Portmanteau Theorem, \eqref{rn_derivative}, and the Cauchy-Schwarz inequality
	$$\PP^*_\mu\left(|Z_t|<\epsilon\right) \le\liminf_n \PP_\mu^{z_n}\left(|Z_t|<\epsilon\right) =
		\liminf_n \E_0^{z_n}\left[\zeta(t) 1_{\{|Z_t|<\epsilon\}}\right]\le$$
		$$\liminf_n\left(\E_0^{z_n}\left[\left(\zeta(t)\right)^2\right]   \PP^{z_n}_0\left(|Z_t|<\epsilon\right)\right)^{1/2}\le \exp\left\{\frac{1}{2}|\mu|^2t\right\}\limsup_n  \left(\PP^{z_n}_0\left(|Z_t|\le\epsilon\right)\right)^{1/2}\le$$
		$$\exp\left\{\frac{1}{2}|\mu|^2t\right\}  \left(\PP^{z}_0\left(|Z_t|\le\epsilon\right)\right)^{1/2}.$$
		The last inequality follows from the Feller property for the drift-less case (\cite{varadhan1985brownian}, Theorem 3.13). We now let $\epsilon\downarrow 0$, and then from property 3 in the definition of the Submartingale Problem follows that 
			$$\PP^*_\mu\left(Z_t=0\right)=0,$$
			for Lebesgue-almost every $t\ge 0$. Identity \eqref{goal} follows.

	\end{proof}
	
	\begin{proof}[Proof of Theorem \ref{Strong Markov for Submartingale Problem}] This follows in the standard way from the uniqueness of the solution to the submartingale problem, using the regular conditional probability measures  for ${\mathcal M}$ given ${\mathcal M}_\tau$ under $\PP^z_\mu$. In particular, Lemma 3.1, and Corollary 3.3 in \cite{varadhan1985brownian} remain true in the presence of a drift, without changing a single word in their proofs. Then the strong Markov property follows, again exactly the same way as in \cite{varadhan1985brownian}, Theorem 3.14.

\end{proof}	

	
	\section{Proof of Theorem \ref{C-hat Feller}}\label{section_C_hat}
	
	\begin{proof}[Proof of Theorem \ref{C-hat Feller}]
		It is sufficient to show  that the $\hat C(S)$-Feller property holds in the case when the drift is zero. Indeed, suppose that for the driftless case the $\hat C(S)$-Feller property holds. By Theorem 1.10 in Bottcher, Schilling, and Wang \cite{bottcher2013levy}, the $\hat C(S)$-Feller property holds for $\{\PP_\mu^z, \forall z\in S\}$ if and only if there exists an increasing sequence of bounded sets $B_n\in{\cal B}(S)$ with $\cup_{n\ge 1}B_n=S$ such that for every $t>0$ and $n\ge 1$
		$$\lim_{|z|\to\infty}\PP_\mu^z(Z(t)\in B_n)=0.$$
		Since we already know that uniqueness holds for the submartingale problem with drift $\mu$, we may assume that $\{\PP_\mu^z,z\in S\}$ is exactly the family we created in the existence part of this paper using Girsanov's theorem. By Theorem 1.10 in \cite{bottcher2013levy}, there exists an increasing sequence of bounded sets $B_n\in{\cal B}(S)$ with $\cup_{n\ge 1}B_n=S$ such that for every $t>0$ and $n\ge 1$,
		$$\lim_{|z|\to\infty}\PP^z_0(Z(t)\in B_n)=0,$$
		where $\{\PP^z, z\in S\}$ is the solution of the submartingale problem without drift. By formula \eqref{rn_derivative}  we have that
		$$\PP^z_\mu(Z(t)\in B_n)=\E^z_0\left[\zeta(t)\mathds{1}_{\{Z(t)\in B_n\}}\right]\le \left(\E^z_0\left[\left(\zeta(t)\right)^2\right]\PP^z(Z(t)\in B_n)\right)^{1/2}$$
		$$=\exp\left\{{\mu^2\over 2}t\right\}\left(\PP_0^z\left(Z(t)\in B_n\right)\right)^{1/2},$$
		which shows that the $\hat C(S)$-Feller property holds for $\{\PP_\mu^z, \forall z\in S\}$. 
		
		In the rest of this proof we shall show that in the case of $\mu=0$ the $\hat C(S)$-Feller property holds. Let $X$ be the process identified in Proposition \ref{Theorem 2.4 in lakner2019roughness}. It has been shown in Williams and Varadhan \cite{varadhan1985brownian} that the Feller property holds, hence the $C_b(S)$-Feller property holds, so by Theorem 1.10 in \cite{bottcher2013levy}, it is sufficient to show that for every $t>0$
		$$\PP_0^z(Z(t)\in B_n)\to 0$$
		as $|z|\to\infty$, where $B_n=\{z\in S: |z|\le n\}$, for $n\ge 1$. Note that
		$$ \PP_0^z(Z(t)\in B_n) = \PP_0^z(Z(t)\in B_n,\tau\le t) + \PP_0^z(Z(t)\in B_n, \tau>t),$$
		where 
		$$\tau=\inf \{t\ge 0: Z(t)\in\partial S\}.$$
		We treat the second term first. It is bounded above by
		$$\PP^z_0(X(t)\in B_n) \le \PP_0^z(|X(t)-z|\ge |z|-n) = { \Pr }(|w(t)|>|z|-n)\to 0,$$
		as $|z|\to\infty$,
		where $w$ is a standard 2-dimensional Brownian motion started at zero.
		Next we treat the first term. Let $T_n=\inf\{t\ge 0: Z(t)\in B_n\}$. Then the first term can be written as a sum of three terms
		$$\PP^z_0\left(Z(t)\in B_n,\tau\le t, |Z(\tau)|<{|z|\over 2}\right) + $$
		$$\PP^z_0\left(Z(t)\in B_n,\tau\le t, T_n\le\tau, |Z(\tau)|\ge{|z|\over 2}\right) + $$
		$$\PP_0^z\left(Z(t)\in B_n,\tau\le t, T_n>\tau, |Z(\tau)|\ge{|z|\over 2}\right)  $$

		Again, we treat the three terms separately. By Proposition \ref{Theorem 2.4 in lakner2019roughness} we have that $Z(\tau)=X(\tau)$, where $X$ is a standard Brownian motion with zero drift started at $z$ under $\PP^z_0$, thus first term is bounded above by
		$$\PP^z_0\left(\tau\le t, |X(\tau)|<{|z|\over 2}\right)\le \PP^z_0\left(\tau\le t, |X(\tau)-z|>{|z|\over 2}\right)\le {\Pr}\left(\max_{s\le t}|w(s)|>{|z|\over 2}\right)\to 0,$$
		as $|z|\to \infty$. The second term is bounded above by
		$$\PP^z_0(T_n\le\tau\le t)\le \PP^z_0(X-z\  \hbox{reaches}\ B_n-z\  \hbox{by time}\ t)= \Pr (w\ \hbox{reaches}\ B_n-z\ \hbox{by time}\ t)\to 0$$
		as $|z| \to\infty$.   Here $w$ is a standard 2-dimensional Brownian motion started at the origin.

		For analyzing the third term we define the stopping time $T_n^\tau=\inf\{t\ge \tau: Z(t)\in B_n\}$. The third term is bounded above by
		$$\PP^z_0\left(\tau<T_n\le t, |Z(\tau)|>{|z|\over 2}\right)\le \PP^z_0\left(\tau<T_n\le \tau+t, |Z(\tau)|>{|z|\over 2}\right)\le$$
		$$ \PP^z_0\left(T_n^\tau\le t, |Z(\tau)|>{|z|\over 2}\right).$$
		By the strong Markov property this can be written as
		$$\int_{\partial S\cap B_{|z|/2}^c}\PP_0^x(T_n\le t)P^z_0(Z(\tau)\in dx).$$
		By the scaling property (Lemma 2.1 in \cite{williams1987local}), the process $\{Z(t), t\ge 0\}$ under $\PP^x_0$ induces the same measure on ${\cal M}$ as $\{|x|Z(t/|x|^2), t\ge 0\}$ induces under $\PP^{x/|x|}_0$, for every non-zero  $x\in S$. 
		Then the above expression can be written as 
		$$\int_{\partial S\cap B_{|z|/2}^c}\PP_0^{x/|x|}\left(|x|^2 T_{n/|x|}\le t\right)\PP_0^z(Z(\tau)\in dx)=$$
		$$\int_{\partial S_1\cap B_{|z|/2}^c}\PP_0^{u_1}\left(|x|^2T_{n/|x|}\le t\right)\PP^z_0(Z(\tau)\in dx) +$$
		$$ \int_{\partial S_2\cap B_{|z|/2}^c}\PP_0^{u_2}\left(|x|^2 T_{n/|x|}\le t\right)\PP_0^z(Z(\tau)\in dx),$$
		where $u_1$ and $u_2$ are the unit vectors $u_1=(1,0)$, and $u_2=(\cos\xi,\sin\xi)$.
		By symmetry it is sufficient to show that the first term converges to 0 as $|z|\to\infty$. If $|z|/2>2n$, then it is bounded above by
		$$\sup_{|x|>|z|/2} \PP_0^{u_1}(|x|^2T_{1/2}\le t)=\PP_0^{u_1}\left({|z|^2\over 4}T_{1/2}\le t\right)\to 0,$$
		as $|z|\to\infty$, which completes the proof of the proposition.
		
	\end{proof}

	\section{Proof of  Theorem \ref{Existence and Uniqueness of Absorbed Process}}

	
	
	
	Before proving the existence part of Theorem \ref{Existence and Uniqueness of Absorbed Process}, we must first establish some preliminary results. Let $\{B_t, t\ge0\}$ be the coordinate mapping process on $C(\R_+,\R^2)$, whose natural filtration is given by $\calW_t = \sigma(B_s, 0\le s \le t)$ for $t \geq 0$, and let $\calW = \sigma(B_s, s\ge 0)$. 
	Recall that $v_i$ is the reflection direction on $\partial S_i$  for $i=1,2$, and let $R$ be the $2\times2$ matrix defined by $R_{ij} = \text{the } i\text{-th component of } v_j$. The following result is adapted from Theorem 3.1 in \cite{williams1983brownian}. 
	

	\begin{prop}\label{Skorokhod problem before hitting the origin} For any $w \in C(\R_+,\R^2)$ with $w(0) \in S$, there exists a unique triple $(\phi,\eta,T_0)$, where $\phi \in C_S$, $\eta \in C(\R_+, [0, +\infty]^2)$ and $T_0: C(\R_+,\R^2) \ra [0, +\infty]$, satisfying the following four conditions,
		\begin{itemize}
			\item[1.] $\phi(t) = w(t) + R\eta(t)$ for each $t \in [0,T_0)$;
			\item[2.] $\phi(t) \neq 0$ for all $t < T_0$ and $\phi(t) = 0$ for all $t \ge T_0$;
			\item[3.] For $j = 1,2$, $\eta(0) = 0$ and $\eta_j(\cdot)$ is non-decreasing and finite for $t \in [0,T_0)$;
			\item[4.] For $j = 1,2$, $\eta_j$ only increases when $\phi(t)$ is on $\partial S_j\backslash\{0\}$.
		\end{itemize}
		Furthermore, we have the following two properties 
		\begin{itemize}
			\item[(i)] $T_0$ is a stopping time on $(C(\R_+,\R^2), \calW, \calW_t)$;
			\item[(ii)] Define the map $\Gamma: C({\mathbb R}_+,{\mathbb R}^2)\mapsto C_S$ such that $\Gamma(w)= \phi$. Then, $\tau_0 \circ \Gamma = T_0$, the map $\Gamma_t  \equiv \Gamma(\cdot)(t)$ is $\calW_{t}$-measurable and $\Gamma$ is $\calW/ \calM$-measurable.
		\end{itemize} 
	\end{prop}
	
	Now for each $z \in S$, let $\hat{\PP}^z_0$ be the unique measure on $(C(\R_+,\R^2), \calW, \calW_t)$ under which $\{B_t, t\ge0\}$ is a standard Brownian motion started at $z$ (the subscript 0 indicates that the Brownian motion has zero drift under $\hat \PP^z_0$). Next, for each $\mu \in \mathbb{R}^2$ and $T \geq 0$, define the measure $\hat {\PP}_{\mu,T}^z$ on $\calW_T$ by
	\begin{equation}
		\frac{d \hat {\PP}_{\mu,T}^z}{d\hat\PP^z_0} = \exp\bigg(\mu(B_T-z)-\frac{1}{2}||\mu||^2T\bigg) ,\nonumber
	\end{equation}
	and define also
	$$\hat B_t=B_t-\mu t,$$
	and note that $\hat B$ is a standard 2 dimensional Brownian motion started at $z$ under {$\hat \PP^z_{\mu , T}$} on $[0,T]$. 
	Then, by Theorem 4.2 in \cite{parthasarathy2005probability} there exists a measure $\hat\PP_\mu^z$ on $\calW$ which coincides with $\hat\PP_{\mu,T}^z$ on $\calW_T$ for all $T \ge 0$. For every $z\in S$ and $\mu\in{\mathbb R}^2$ the process $B$ is a Brownian motion with drift $\mu$ started at $z$ under the probability measure $\hat \PP^z_\mu$.
	Moreover, since by Proposition \ref{Skorokhod problem before hitting the origin}, $\Gamma$ is a measurable map from $(C(\R_+, \R^2), \calW, \calW_{t})$ to $(C_S, \calM, \calM_{t})$, we may denote by $\bPm^z$ the measure induced on $\cal M$ by the mapping $\Gamma$ under $\hPm^z$, i.e.,
	\begin{equation}\label{construction of the measure of the absorbed process}
		\bPm^z(A) \equiv \hPm^z(\Gamma^{-1}(A)), \text{ for each } A \in \calM. 
	\end{equation}
	
	We next prove that the family of measures $\{\bPm^z, z \in S\}$ is a solution to the absorbed process problem of Definition \ref{definition of absorbed process}. Since conditions $1$ and $3$ of Definition \ref{definition of absorbed process}  are trivially satisfied by $\{\bPm^z, z \in S\}$, it remains to prove that condition $2$ is satisfied as well. This is achieved in Lemma \ref{satisfying the submartingale property}.
	
	\begin{lem}\label{satisfying the submartingale property}
		Let the family of measures $\{\bPm^z,  z \in S\}$ be defined as in (\ref{construction of the measure of the absorbed process}) and let $Z$ be the coordinate-mapping process on $(C_S, \calM,\calM_{t})$. 
		Then, the process 
		\begin{equation}
		\bigg\{f(Z({t\wedge \tau_0})) - \int_0^{t\wedge \tau_0} \mu \cdot \nabla f(Z(s))ds -\frac{1}{2} \int_0^{t\wedge \tau_0} \Delta f(Z(s))ds, t \ge 0\bigg\}
		\label{submtg}
		\end{equation}
		is a submartingale on $(C_S, \calM, \calM_{t}, \bPm^z)$, for each $f \in C_b^2(S)$ such that $D_i f \ge 0$ on $\partial S_i$ for $i = 1, 2$.
	\end{lem}
	\begin{proof}
		For each $w \in C(\R_+, \R^2)$, let  $\phi(w)=\Gamma(w)$ and note that by Proposition \ref{Skorokhod problem before hitting the origin} we may write $\phi(t) = w(t) + R\eta(t)$ for all $t \ge 0$. Now{, on $(C(\R_+, \R^2), \calW_t, \calW)$, we consider the process $\{\phi(t), t\ge 0\}$ with $\phi(t) = \Gamma_t(w)$ for any $w \in C(\R_+, \R^2)$.} {Recall that the coordinate mapping process $\{B_t, t \ge 0\}$ is a Brownian motion on $(C(\R_+, \R^2), \calW_t, \calW)$ under $\hat{\PP}^z$, and $\{B_t - \mu t, t \ge 0\}$ is a Brownian motion on $(C(\R_+, \R^2), \calW_t, \calW)$ under $\hPm^z$. Notice, by Theorem 1 in \cite{williams1985reflected}, that $R\eta(t)$ is of finite variation on $[0, t\wedge T_0]$ for any $t\ge0$, and $w(t) = B(w)(t)$, we get that $\{\phi(t), t\ge 0\}$ is a semimartingale under $\hat{\PP}^z$. On the other hand, the Girsanov transform keeps the semimartingale property, so $\{\phi(t), t\ge 0\}$ is also a semimartingale under $\hPm^z$.} Hence, for each $f \in C_b^2(S)$ such that $D_i f \ge 0$ on $\partial S_i$ for $i = 1, 2$, {we use It\^{o}'s formula under $\hPm^z$} and get
		\begin{eqnarray}
			f(\phi({t\wedge T_0})) -f(\phi(0))&=& \sum_{i=1}^2 \int_0^{t\wedge T_0} {\partial f \over \partial x_i} (\phi(s))d(w_i(s)-\mu_i s) + \int_0^{t\wedge T_0} \mu \cdot \nabla f(\phi(s))ds  \nonumber \\
			&+& \int_0^{t\wedge T_0} (D_1 f(\phi(s)), D_2 f(\phi(s))) \cdot d\eta(s)+ \frac{1}{2} \int_0^{t\wedge T_0} \Delta f(\phi(s))ds \nonumber
		\end{eqnarray}
		Since we have by Proposition \ref{Skorokhod problem before hitting the origin} that for $i=1,2$, 
		$$
		d\eta_i(s) = \mathds{1}_{\{\phi(s) \in \partial S_i\backslash \{0\}\}}d\eta_i(s),~~s \geq 0,
		$$
		and by the assumption on $f$ that for $i=1,2,$
		$$
		D_i f(\phi(s))  \mathds{1}_{\{\phi(s) \in \partial S_i \backslash \{0\}\}} \ge 0,~~s \geq 0,
		$$
		it follows that the process
		$$
		\left\{ \int_0^{t\wedge T_0} (D_1 f(\phi(s)), D_2 f(\phi(s))) \cdot d\eta(s),  t\ge 0 \right\}
		$$
		is increasing. On the other hand, since { $\{B_t - \mu t, t \ge 0\}$ is a Brownian motion under $\hPm^z$}, the process
		$$
		\left\{ \sum_{i=1}^2 \int_0^{t\wedge T_0} {\partial f \over \partial x_i} (\phi(s))d(w_i(s)-\mu_i s), t\ge 0 \right\}
		$$
		is a martingale under $\hat{\PP}_\mu^z$, so 
		\begin{eqnarray}
			&&f(\phi({t\wedge T_0})) - \int_0^{t\wedge T_0} \mu \cdot \nabla f(\phi(s))ds -  \frac{1}{2} \int_0^{t\wedge T_0} \Delta f(\phi(s))ds \nonumber \\
			&=& f(\phi(0)) + \sum_{i=1}^2 \int_0^{t\wedge T_0} {\partial f \over \partial x_i} (\phi(s))d(w_i(s)-\mu_i s) \nonumber \\
			&+& \int_0^{t\wedge T_0} (D_1 f(\phi(s)), D_2 f(\phi(s))) \cdot d\eta(s) \nonumber
		\end{eqnarray}
		is a submartingale under $\hat{\PP}_\mu^z$. It  follows from \eqref{construction of the measure of the absorbed process} that the process under \eqref{submtg}  is also a submartingale under the {induced} measure $\bar{\PP}_\mu^z$.
	\end{proof}

	\begin{proof}[Proof of the existence part of Theorem \ref{Existence and Uniqueness of Absorbed Process}] The existence of a solution to the absorbed process problem  follows from  Lemma \ref{satisfying the submartingale property}.
	\end{proof}

	\noindent{\it Proof of the uniqueness part of Theorem \ref {Existence and Uniqueness of Absorbed Process}  }. The proof of the uniqueness of the solution to the absorbed process problem is 
very similar to that of the solution of the submartingale problem, hence in this section we shall state the appropriate lemmas and indicate the necessary changes in order to adapt the proofs in Section \ref{Existence} to the absorbed process problem.
	
	\begin{lem}\label{Ocuppation time on the boundary is zero again}
		Let $ \alpha \in{\mathbb R}$ arbitrary, and suppose that $\{\PP_\mu^{z,0}, z\in S\}$ is a solution to the absorbed process problem with drift $\mu \in \R^2$. Then, for all $z \in S$, 
		
		$$\E_\mu^{z,0}\bigg[\displaystyle \int_0^{\tau_0} \mathds{1}_{\{Z(t) \in \partial S\}} dt \bigg] = 0.$$
		
	\end{lem}
	
	\begin{proof} The proof is almost identical to that of Lemma \ref{Occupation:time} with the modification that all processes must be stopped at $\tau_0$.

	\end{proof}

	\begin{lem}\label{Existence of BM under drift measure again}
		Let $ \alpha \in{\mathbb R}$ be arbitrary. Suppose that $\{\PP_\mu^{z,0}, z\in S\}$ is a solution to the absorbed process  problem with drift $\mu \in \R^2$. Then there exists a process $X$ on $(C_S, \calM, \calM_t)$ such that for all $z \in S$, $X$ is a Brownian Motion with drift $\mu$  under $\PP_\mu^{z,0}$ started at $z$ and stopped at $\tau_0$. In addition, $Y = Z-X$ is flat on $[\sigma_n^\delta, \tau_n^\delta]$.
	\end{lem}
	
	\begin{proof}
		The proof is very similar to the proofs  of Lemmas \ref{w:delta}, \ref {construct the Brownian motion using the convergence}, and   \ref{Existence of BM under drift measure}, with the difference that all processes must be stopped at $\tau_0$.
	\end{proof}
	
	\noindent{\it Proof of the uniqueness part of Theorem \ref{Existence and Uniqueness of Absorbed Process}.} The proof is basically a copy of the proof of the uniqueness of the solution to the submartingale problem (Lemma \ref{unique:main} and the proof of Theorem \ref{Exsistence and Uniqueness of Submartingale Problem}).      The necessary changes in order to adapt that proof to the present situation are the following:
	\begin{itemize}
		\item All processes must be stopped at $\tau_0$;
		\item The constraint that the test function $f$ is constant in a neighborhood of the vertex must be erased (above \eqref{decomposition1} in the adapted proof;
		\item At the end of the proof, instead of using the uniqueness of the solution to the submartingale problem with no drift, we must use the uniqueness of the solution to the absorbed process problem with no drift (\cite{varadhan1985brownian}, Theorem 2.1).\hfill$\square$
	\end{itemize}

	\section{Proof of Theorems \ref{Hitting prob for alpha <= 0}, \ref{Hitting prob for alpha >= 1}, and Proposition \ref{stays}} \label{sec:hitting}

	We need to recall Proposition \ref{Skorokhod problem before hitting the origin}. According to this proposition, for every $w\in C({\mathbb R}_+,{\mathbb R}^2)$ with $w(0)\in S$ there exists a triple
	$(\phi,\eta,T_0)$ such that items 1-4 and (i),(ii) hold. Since $B$ is the coordinate mapping process on  $C({\mathbb R}_+,{\mathbb R}^2)$, we may replace $w$ in item 1 with $B$, and write 
	$$\phi(t)=B(t)+R\eta(t),\quad t\in[0,T_0).$$
	We also know that $B$ is a 2 dimensional Brownian motion with drift $\mu$ started at $z$ under $\hat \PP_\mu^z$ for every $z\in S$, hence we can write
	\begin{equation} \phi(t)=z+ W(t)+\mu t +R\eta(t),\quad t\in[0,T_0),\label{replace}\end{equation}
	where $W$ is a standard  2 dimensional Brownian motion started at zero under $\hat\PP_\mu^z$. 
	The measure $\PP^{z,0}_\mu$ was defined as the measure induced by $\Gamma$ on ${\mathcal  M}$ under $\hat\PP^z_\mu$. From this and from $\tau_0\circ \Gamma=T_0$ follows that
	\begin{equation}\PP_\mu^{z,0}(\tau_0<\infty) = \hat\PP_\mu^z(T_0<\infty).\label{change}
	\end{equation}
	We shall use \eqref{change} repeatedly in the coming proofs.
	
	\begin{proof}[Proof of Theorem \ref{Hitting prob for alpha <= 0}]
		By Theorem 2.2 in \cite{varadhan1985brownian}, we have that $\PP_0^{z,0}(\tau_0 = \infty) = 1$, thus by \eqref{change} also $\hat\PP_0^{z}(T_0 = \infty) = 1$ . For every $n\in{\mathbb N}_+$ the measures $\hat \PP_0^z$ and $\hat \PP^z_\mu$ are mutually absolutely continuous on ${\mathcal W}_n$, so $\hat\PP_0^{z}(T_0 <n) = 0$ implies  $\hat\PP_\mu^{z}(T_0 <n) = 0$. Then $\hat\PP_\mu^{z}(T_0 =\infty) = 1$ follows, and this and \eqref{change} gives the required result.
	\end{proof}

	\begin{proof}[Proof of Theorem \ref{Hitting prob for alpha >= 1}]
		First we are going to show \eqref{finite}.  By the $\alpha\ge 1$ condition  there exists a vector $ b\in \R^2$ such that $b \cdot z < 0$ for all $z\in S, z\neq 0$, and $b\cdot v_i \ge 0$ for $i=1,2$. 
Indeed, if $\alpha\ge 1$ then ${\rm co}(-v_1,-v_2)$ is either a line containing $S$ within one side, or it is a wedge with angle less than $\pi$ containing $S$. In either case 		
		the existence of such a vector follows. Then, by identity \eqref{replace}, for each $z \in S$,
		$$
		0 \ge b\cdot \phi(t) = b\cdot z +b \cdot W_t + b \cdot v_1 \eta_1(t) + b\cdot v_2 \eta_2(t) + b\cdot \mu t \ge b \cdot z + b\cdot W_t +b \cdot \mu t, 
		$$
		for $t < T_0$, $\hat\PP^{z}_{\mu}$-a.s., and  so
		$$
		\hat\PP_\mu^{z}(0 \ge b \cdot z + b\cdot W_t +b \cdot \mu t,\ t < T_0) = 1.
		$$
		Therefore,
		\begin{eqnarray}
			&&\hat\PP_\mu^{z}(0 \ge b \cdot z + b\cdot W_t +b \cdot \mu t,\ t < \infty)\nonumber\\
			&\ge& \hat\PP_\mu^{z}(\{0 \ge b \cdot z + b\cdot W_t +b \cdot \mu t, t < T_0\}\cap \{T_0 = \infty\})\nonumber\\
			&=& \hat\PP_\mu^{z}(T_0 = \infty). \nonumber
		\end{eqnarray}
		This implies 
		\begin{equation}
		\hat\PP_\mu^{z}(-b \cdot z \ge b\cdot W_t +b \cdot \mu t, t < \infty)
		\ge \hat\PP_\mu^{z}(T_0 = \infty).
		\label{ineq}\end{equation}
		However,
		$\hat\PP_\mu^{z}(-b \cdot z \ge b\cdot W_t +b \cdot \mu t,\ t < \infty) <1$, and together with \eqref{change} this proves  the result.
Suppose now that \eqref{condition} also holds, in addition to $\alpha\ge 1$. Then $\rm{co}(v_1,v_2,\mu)$ is either a wedge with angle less than $\pi$, or a half-space, or a line. Then the same is true for  $\rm{co}(-v_1,-v_2,-\mu)$, and if it is a wedge or a half-space then it contains $S$, and if it is a line then it contains $S$ in one side. In all cases we can select $b$ so that in addition to 	$b \cdot z < 0$ for all $z\in S, z\neq 0$, and $b\cdot v_i \ge 0$ for $i=1,2$, we also have $b\cdot\mu\ge 0$. In that case $\hat\PP_\mu^{z}(-b \cdot z \ge b\cdot W_t +b \cdot \mu t,\ t < \infty) =0$, so \eqref{ineq} and \eqref{change} imply \eqref{reaches}.
	\end{proof}

	\begin{proof} [Proof of Proposition \ref{stays}.] Suppose first that $0 < \xi \leq \pi/2$. Let $z \in S^0$ and $\mu \in \mathbb{R}^2$ be given by $\mu = ||\mu||(\cos \eta, \sin \eta) \neq 0$, where  $\eta \in (0,\xi)$. Next, set 
\begin{eqnarray}
X_t &=& z + B_t + \mu t,~t \geq 0,
\end{eqnarray}
where $B_t$ is a standard $2$-d Brownian motion.

Now translate the origin of the coordinate axes to $z$ and then rotate the axis in a counterclockwise direction by the angle $\eta$. By the translational and rotational invariance of Brownian motion, in these new coordinates the process $X$ may be written as
\begin{eqnarray*}
\hat{X}_t = \hat{B}_t + \|\mu \|\hat{e}_1 t, ~t \geq 0,
\end{eqnarray*}
where $\hat{B}$ is a standard Brownian motion and $\hat{e}_1=(1,0)$.  Next, in the new coordinate system denote by $\mathcal{L}_1$ the line corresponding $\partial S_1$, and by $\mathcal{L}_2$ the line corresponding to $\partial S_2$. Then, in the new coordinate system the interior of $S$ may be expressed as
\begin{eqnarray*}
\mathcal{S}^0 &=& \{\hat{z} \in \mathbb{R}^2 :\mathcal{L}_1(\hat{z}_1) < \hat{z}_2 < \mathcal{L}_2(\hat{z}_1) \},
\end{eqnarray*}
where $\mathcal{L}_i(\hat{z}_1)$ is a coordinate uniquely determined by the relation $(\hat z_1, \mathcal{L}_i(\hat{z}_1))\in\mathcal{L}_i$, for $i=1,2$. Hence, in order to complete the proof for the case of $0 < \xi \leq \pi/2$ it suffices to show that
\begin{eqnarray*}
P( \mathcal{L}_1(\hat{X}_t^1) < \hat{X}_t^2 < \mathcal{L}_2(\hat{X}_t^1), t \geq 0 ) &>&0.
\end{eqnarray*}

First note that since $0 < \eta < \xi \leq \pi/2$, we may write
\begin{eqnarray*}
     \mathcal{L}_1(\hat{z}_1) &=&  -a - b \hat{z}_1 ~~~\mathrm{and}~~~ \mathcal{L}_2(\hat{z}_1) ~=~  c + d \hat{z}_1,
\end{eqnarray*}
for $a,b,c,d > 0$. Hence,
\begin{eqnarray}
&& \{\mathcal{L}_1(\hat{X}_t^1) < \hat{X}_t^2 < \mathcal{L}_2(\hat{X}_t^1), t \geq 0\} \label{set:equality}\\
&=& \{ -a - b\| \mu \| t - b \hat{B}_t^1 < \hat{B}_t^2 < c + d \|\mu \|t + d \hat{B}_t^1 , t \geq 0\}. \nonumber
\end{eqnarray}
From \eqref{set:equality}, it now follows after some algebra that 
\begin{eqnarray*}
&&\{\hat{B}^1_t > \max(-a/2b,-c/2d) - (\| \mu\|/2)t , t \geq 0 \}\\
&& \cap  \{ -(a/2) - (b/2)\| \mu \| t  < \hat{B}_t^2 < (c/2) + (d/2) \|\mu \|t  , t \geq 0\} \\
&\subseteq& \{\mathcal{L}_1(\hat{X}_t^1) < \hat{X}_t^2 < \mathcal{L}_2(\hat{X}_t^1), t \geq 0\}.
\end{eqnarray*}
Next, by the independence of $\hat{B}^1$ and $\hat{B}^2$, we have that
\begin{eqnarray*}
&&P( \mathcal{L}_1(\hat{X}_t^1) < \hat{X}_t^2 < \mathcal{L}_2(\hat{X}_t^1), t \geq 0 ) \\
&\geq&P(\hat{B}^1_t > \max(-a/2b,-c/2d) - (\| \mu\|/2)t , t \geq 0 )\\
&& \times P( -(a/2) - (b/2)\| \mu \| t  < \hat{B}_t^2 < (c/2) + (d/2) \|\mu \|t  , t \geq 0).
\end{eqnarray*}
However, since $a,b,c,d,\| \mu \| > 0$, it follows by (4.3) of Doob \cite{doob1949heuristic} that the 2 probabilities on the righthand side above are greater than zero. This completes the proof for the case of $0 < \xi \leq \pi/2$.

Now suppose that $\pi/2 < \xi < \pi$ and $z \in S^0$ and $ 0 < \eta < \xi $. In this case we show that there exists a wedge $\bar{S} \subset S$ such that $z \in \bar{S}^0$ and $P(X_t \in \bar{S}^{0} , t \geq 0 ) > 0$, which is sufficient to complete the proof. First, suppose that $0 < \eta < \pi/2$. In this case the wedge $\bar{S}$ can be defined in the usual way by setting $\xi = \pi/2$ and placing the vertex of $\bar{S}$ at a point $\bar{z} \in S^{0}$ such that $\bar{z} < z$. The results above then yield the desired result. Next, suppose $\pi/2 \leq \eta < \xi$. Then, set $\bar{S}=v+ \{r \geq 0, \pi - (\xi+ \eta)/2 \leq \theta \leq (\xi+\eta)/2  \}$, where the vertex $v=(v_1,v_2) \in S^{0}$ is such that $v_1=z_1$ and $v_2 < z_2$. The results above again yield the desired result.
\end{proof}



	\section{Appendix}
	
	\begin{lem}\label{f:function}
		There exists a function $f_{\varepsilon,C}\in C_b^2(S)$ satisfying
			\begin{equation*} f_{\varepsilon,C}(x,y)= 
			\begin{cases}
				0,\ {\rm if}\ (x,y)\in S\setminus S^{\varepsilon/3},\\
				y,\ {\rm if}\ (x,y)\in S^{2\varepsilon/3}, y\le C,
			\end{cases}\end{equation*}
		such that in addition $f_{\varepsilon,C}(x,0) =0$ for all $x\ge 0$, and $D_i f_{\varepsilon, C}\ge 0$ on $\partial S_i$.
		\end{lem}
	
	\begin{proof}
		Let $h_1\in C_b^2({\mathcal R})$ such that $h_1(x) \ge 0$ for all $x\in{\mathbb R}$
		and 
			\begin{equation*} h_1(x)= 
			\begin{cases}
				0,\ {\rm if}\ x\le\varepsilon/3,\\
				1,\ {\rm if}\ x\ge 2\varepsilon/3.
			\end{cases}\end{equation*}
		Let  $h_2\in C_b^2({\mathcal R})$ such that $h_2(y)=y$ if $y\le C$.
		Then
		$$f_{\varepsilon,C} =h_1\left(x -\frac{y}{\tan\xi}\right)h_2(y)\quad (x,y)\in S$$
		satisfies the requirements of the lemma.
		 Note that for any $\delta>0$ we have $(x,y)\in S^\delta$ if and only if $x-y/\tan\xi\ge\delta.$ Using this fact repeatedly one can verify the above statement by straightforward calculation.
	\end{proof}
	
	\begin{rmk}\label{aug} There are slightly different definitions for the term  ``augmented filtration" in the literature; we use this term as defined in \cite{RogersWilliams}, Definition II.67.3.
 Let $\PP$ be an arbitrary probability measure on ${\mathcal M}$, and let $(C_S, {\mathcal F}, ({\mathcal F}_t), \PP)$ be the  augmentation of the probability space $(C_S,{\mathcal M},({\mathcal M}_t),\PP)$, in the above sense.
It is known that right-continuous martingales (submartingales) on $(C_S,{\cal M},({\cal M}_t),\PP)$ are also right-continuous martingales (submartingales) on $(C_S,{\cal F},({\cal F}_t),\PP)$ (Lemma II.67.10 in \cite{RogersWilliams}). The propability measure $\PP$ in the second probability space is the extention of $\PP$ from ${\mathcal M}$ to ${\mathcal F}$, without changing the notation.  Also, it follows from the martingale characterization of Brownian motion (Theorem 3.3.16 in \cite{karatzas2012brownian}) that a Brownian motion on $(C_S,{\cal M},({\cal M}_t),\PP)$ is also a Brownian motion on $(C_S,{\cal F},({\cal F}_t),\PP)$.
\end{rmk}

	\bibliographystyle{plain}
	\bibliography{references}
\end{document}